\documentclass{amsart}
\usepackage{amssymb}
\usepackage{amsthm}
\usepackage{graphics}
\newtheorem*{reftheorem}{Theorem}
\newtheorem*{openprob}{Open Problem}
\newtheorem{introtheorem}{Theorem}
\newtheorem{theorem}{Theorem}
\newtheorem{lemma}{Lemma}[section]
\newtheorem{construction}{Construction}

\newtheorem{introprop}{Proposition}
\newtheorem{prop}{Proposition}[section]

\numberwithin{equation}{section}
\title{The width-volume inequality}
\author{Larry Guth}
\address{Department of Mathematics, Stanford, Stanford CA, 94305 USA}
\email{lguth@math.stanford.edu}

\begin{document}
\begin{abstract}
We prove that a bounded open set $U$ in $\mathbb{R}^n$ has k-width
less than $C(n) \textrm{ Volume}(U)^{k/n}$.  Using this estimate,
we give lower bounds for the k-dilation of degree 1 maps between
certain domains in $\mathbb{R}^n$.  In particular, we estimate
the smallest (n-1)-dilation of any degree 1 map between two
n-dimensional rectangles.  For any pair of rectangles, our
estimate is accurate up to a dimensional constant $C(n)$.  We
give examples in which the (n-1)-dilation of the linear map is
bigger than the optimal value by an arbitrarily large factor.
\end{abstract}

\maketitle

This paper proves some estimates having to do with the areas of
k-dimensional surfaces in Euclidean space.  We deal with two
problems.  First, suppose that $U$ is a bounded open set in
$\mathbb{R}^n$.  We consider the problem of sweeping out $U$ with
k-dimensional surfaces, trying to arrange that the volumes of all
the surfaces are as small as possible.  Depending on the geometry
of $U$, we give upper and lower bounds for the possible volumes
of the surfaces.  In particular, we construct a family of
k-dimensional surfaces sweeping out $U$ so that each surface has
volume bounded by $C(n) \textrm{Volume}(U)^{k/n}$.  The next
question concerns mappings from one open set to another - for
example from the unit cube to a long thin cylinder.  After we fix
a domain and a range, we consider the problem of finding a degree
1 mapping which stretches the volumes of k-dimensional surfaces
as little as possible.  For certain pairs of (n-dimensional)
rectangles, we show that the linear mapping stretches the
k-dimensional surfaces far more than necessary.  We give upper
and lower bounds for the minimal amount of stretching by any
degree 1 map.  When $k = n-1$, these upper and lower bounds match
up to a constant factor.

\vskip5pt
{\bf The definition of k-width}

As a first approximation to the definition of width, we define a
linear version of k-dimensional width. Let $U$ be a bounded open
set in $\mathbb{R}^n$.  For each (n-k)-plane $P$ through the
origin, let $F(P)$ denote the family of all k-planes perpendicular
to $P$.  Define the width of $F(P)$ to be the maximum volume of the
intersection of $U$ with any of the k-planes in $F(P)$.  Then define
the linear k-width of $U$ to be the minimum width of $F(P)$ as $P$
varies among all the (n-k)-planes through the origin.

The width considered in this paper is a non-linear generalization
of the definition above.  Instead of families of parallel
k-planes, we consider families of k-dimensional surfaces.  The
surfaces we consider will be oriented relative k-cycles in
$U$.  For the reader not familiar with k-cycles, there is no harm
in picturing k-dimensional submanifolds of $U$ with boundary in
$\partial U$.  We define the k-width of a family $F$ to be the
largest k-volume of any of the k-cycles in $F$.  In order to define the
k-width of $U$, we look at families of cycles that ``sweep out'' $U$. 
Morally, a closed (n-k)-dimensional family $F$ of k-cycles can be
glued together to form a single n-dimensional cycle.  If this
n-dimensional cycle has a non-zero homology class, then we say
that $F$ ``sweeps out'' $U$.  For example, if $\pi$ is a generic
PL map from $U$ to $\mathbb{R}^{n-k}$, then the fibers
$\pi^{-1}(y)$ form a family of k-cycles sweeping out $U$,
parametrized by $y \in \mathbb{R}^{n-k}$.  We define the k-width
of $U$ to be the smallest k-width of any family of k-cycles
sweeping out $U$.  Because the definition doesn't involve planes,
it also makes sense if we replace $U$ by any compact oriented
Riemannian manifold.

Mathematicians working on geometric measure theory began to look
at families of cycles in the 1960's.  In the unpublished paper
\cite{A}, Almgren used such families as a tool to construct
minimal cycles on a Riemannian manifold using Morse-theoretic
arguments.  A good reference for this material is the first
chapter of Pitts's book \cite{P}.  Gromov had the idea to use
families of cycles as a way of describing the size of a
Riemannian manifold $(M,g)$.  He sketched his ideas about this
subject in section F of appendix 1 of his long paper on metric
geometry \cite{G1}.  In this section he essentially gave the
definition above. 

Because the space of all k-cycles is infinite dimensional, it
takes some work to prove that the k-width of an open set is not
zero.  The first proof of this fact is essentially due to
Almgren.  Gromov pointed out that Almgren's work establishes the
exact k-width of the unit n-sphere: the k-width of the unit
n-sphere is equal to the volume of the unit k-sphere.  Almgren's
proof requires a substantial amount of geometric measure theory. 
Gromov also gave an elementary lower bound for the k-width of the
unit n-sphere.  Using Gromov's proof, it's not hard to estimate
the k-width of simple shapes like rectangles.

\begin{introprop} Let $R$ be an n-dimensional rectangle with
dimensions $R_1 \le ... \le R_n$.  Then the k-width of $R$ is
at least $c(n) R_1 ... R_k$ and at most $R_1 ... R_k$.
\end{introprop}

(It seems reasonable to guess that the k-width of a rectangle is
$R_1 ... R_k$, but the exact value of the k-width is unknown.)

\vskip5pt

{\bf The width-volume inequality}

The first theorem of this paper is an upper bound for the
k-width of sets with small volume in Euclidean space.

\begin{introtheorem} (The width-volume inequality) Let $U$ be a
bounded open set in $\mathbb{R}^n$ with volume $V(U)$ and k-width
$W_k(U)$.  Then $W_k(U) < C(n) V(U)^{k/n}$.
\end{introtheorem}

To prove the theorem, we have to construct a family of cycles
that sweep out $U$ in an efficient way.  The first approach one
might try is to use parallel planes at a well-chosen angle, as in
the definition of linear k-width.  This approach can fail because
of the Kakeya phenomenon.  As proven by Besicovitch, there are
open sets $U$ in $\mathbb{R}^2$ with arbitrarily small area
containing a unit line segment in every direction.  These sets
have linear 1-width at least 1 and arbitrarily small area.  A
good reference for Besicovitch sets is Wolff's article \cite{Wol}. 
While I was revising this paper, I learned that taking parallel
k-planes does work when $k > n/2$.  This result was proven by
Falconer in the interesting short paper \cite{Fal} (in slightly
different language).  In the last section, we briefly explain
Falconer's proof, which is based on Fourier analysis.  For the
intermediate range $2 \le k \le n/2$, I don't know if the linear
k-width of a set $U$ can be bounded in terms of its volume. 
There are some more comments in the open problem section.

We now sketch the proof of our theorem, which deals with all
values of k.  Because of a scaling argument, we can assume that
the volume of $U$ is 1.  The first step of the argument is to find
a lot of k-planes that meet $U$ in a small volume.  We find these
planes by an averaging trick.  Let $S_0$ be the k-skeleton of the
unit lattice in $\mathbb{R}^n$.  We consider the translations of
$S_0$ by a vector $x \in [0,1]^n$.  On average, the translation
of $S_0$ meets $U$ in a region of volume $n \choose k$.  Therefore,
we can choose a translate $S$ of $S_0$ which meets $U$ in a
controlled volume.  We can then control the k-volume of any
k-cycle lying in the skeleton $S$.

It is not possible to sweep out $U$ with cycles lying in the
skeleton $S$, because any family of cycles sweeping out $U$ must pass
through every point of $U$.  But it turns out to be possible to
sweep out $U$ by a family of k-cycles each of which lies in $S$ except
for a subset of controlled volume.

\begin{introprop} For any bounded open set $V \in \mathbb{R}^n$
(of any volume), there is a family of k-cycles sweeping out $V$ so
that each cycle lies in $S$, except for a subset of volume at
most $C(n)$.
\end{introprop}

We include some pictures to indicate how such a family might look
for $k=1$, $n=2$.  The thin lines denote the 1-skeleton $S$ and the
thick lines denote a 1-cycle in our family.

\vskip5pt

\includegraphics{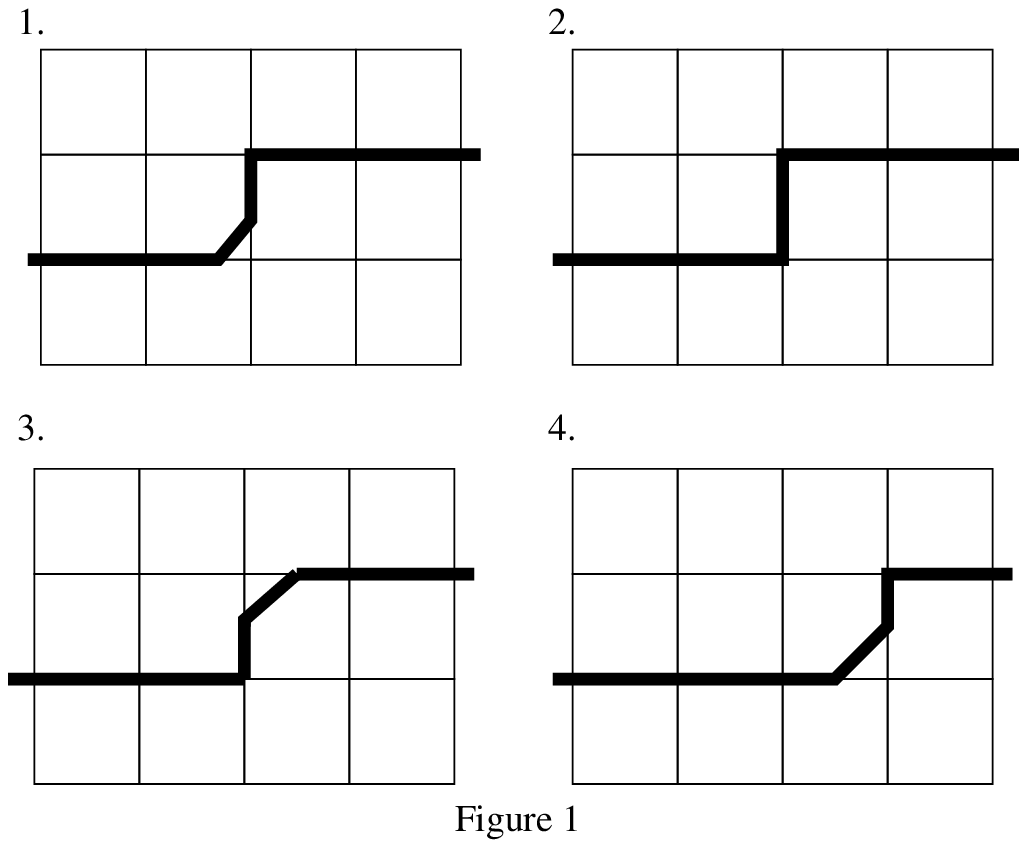}

The general case is somewhat harder than the case $k=1, n=2$.  In
general, the family of k-cycles is constructed by starting with a
family of parallel k-planes transverse to $S$, and then ``bending''
them so that almost all of the volume of each k-plane is pushed
into the skeleton $S$.  We call this construction ``bending planes
around a skeleton''.

\vskip5pt
 
{\bf Area-contracting maps between rectangles}

In the second half of the paper, we apply the width-volume
inequality to estimate the k-dilations of degree 1 maps.  Recall
that k-dilation measures how much a mapping stretches
k-dimensional volumes.  If a map $f$ takes any k-dimensional
manifold with volume $V$ to an image with volume at most $\lambda V$,
then $f$ has k-dilation at most $\lambda$.

The second main problem of this paper is to estimate the infimal
k-dilation of all degree 1 maps from a rectangle $R$ to another
rectangle $S$.  This innocuous-sounding problem has turned out to
be much more complicated than I expected.  When I first
approached the problem, I guessed that a linear diffeomorphism
from $R$ to $S$ would give at least roughly the smallest
k-dilation.  My guess was wrong.  Let us define $D_k(R,S)$ to be
the infimal k-dilation of any degree 1 map from $R$ to $S$
(taking the boundary of $R$ to the boundary of $S$).  For comparison,
let us define $Lin_k(R,S)$ to be the smallest k-dilation of a
linear diffeomorphism from $R$ to $S$.

\begin{introprop} For each $n \ge 3$ and each k in the range $2
\le k \le n-1$, there are pairs of n-dimensional rectangles
$(R,S)$ so that the ratio $Lin_k(R,S)/D_k(R,S)$ is arbitrarily large.
\end{introprop}

For example, if the rectangle $R$ has dimensions $\epsilon \times 1
\times 1$, and the rectangle $S$ has dimensions $\epsilon \times
\epsilon \times \epsilon^{-1}$, then $Lin_2(R,S) =
\epsilon^{-1}$.  On the other hand, there is a non-linear degree
1 map from $R$ to $S$ with 2-dilation less than 1000, regardless of
$\epsilon$.  I call this map the snake map because it somewhat
resembles a snake uncoiling.

We take a little time to describe this map.  The snake map does
not have any analogue in 2 dimensions, but there is a map related
to it.  Let $U$ be the unit square, and let $A
\subset U$ be the shape in Figure 2.

\vskip5pt

\includegraphics{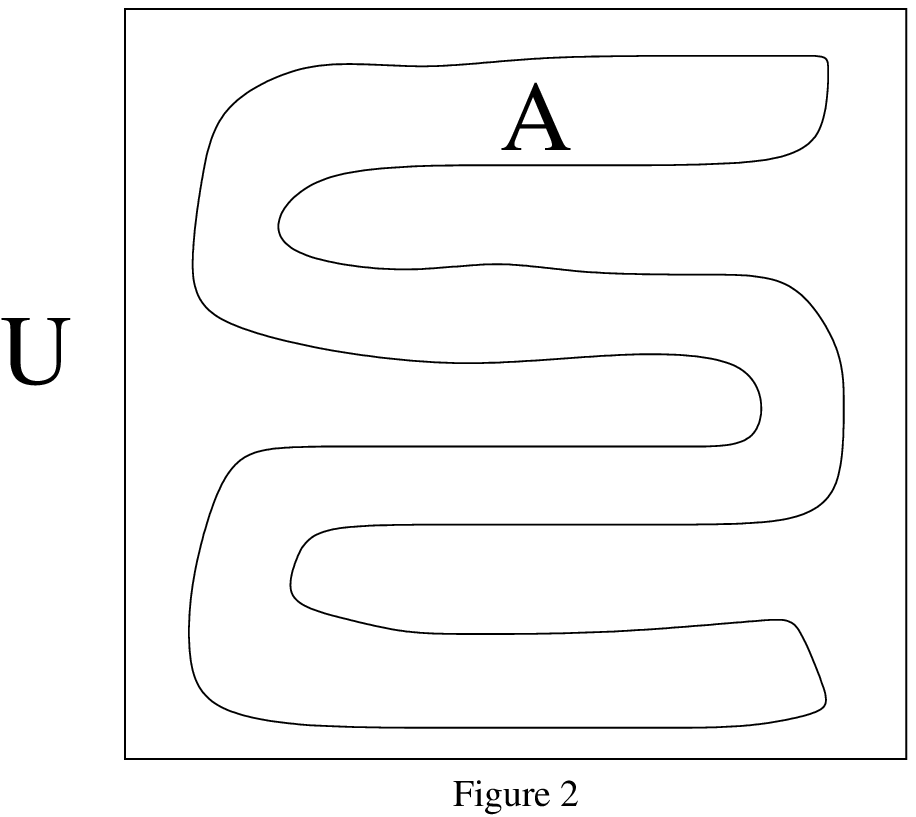}

The set $A$ is bilipschitz to the rectangle $[0, \epsilon]
\times [0, \epsilon^{-1}]$, and it snakes back and forth across $U$
roughly $\epsilon^{-1}$ times.  Let $A^c$ denote the complement
of $A$ in $U$.  The first map that we consider is a retraction $\phi$
of $U$ onto $A$, which maps $A^c$ onto $\partial A$.  The 1-dilation
of $\phi$ is roughly $\epsilon^{-1}$ and the 2-dilation of $\phi$
is exactly 1.

We now turn to three dimensions.  The rectangle $R$ is equal to
$[0, \epsilon] \times U$ and the rectangle $S$ is bilipschitz to
$[0, \epsilon] \times A$.  We can get a degree 1 map from $R$ to $S$
by first retracting $R$ onto $[0, \epsilon] \times A$ and then
using the bilipschitz equivalence of $[0, \epsilon] \times A$ with $S$.

The most obvious retraction from $R$ onto $[0, \epsilon] \times
A$ is $id \times \phi$, where $id$ denotes the identity map from
$[0, \epsilon]$ to itself.  This retraction has 2-dilation
roughly $\epsilon^{-1}$.  Using this retraction, we get a degree
1 map from $R$ to $S$ with 2-dilation roughly $\epsilon^{-1}$,
slightly larger than the 2-dilation of the linear map.

The trick in the construction of the snake map is to improve the
retraction from $R$ to $[0, \epsilon] \times A$.  The improved
retraction takes place in two steps.  We first retract $R$ onto the
union $(\{ 0 \} \times U) \cup ([0, \epsilon] \times A)$.  We
then retract this set onto $[0, \epsilon] \times A$.  The set
$(\{ 0 \} \times U) \cup ([0, \epsilon] \times A)$ resembles a
snake sitting on a piece of cardboard.  The first retraction can
be done with 1-dilation roughly 1, and hence 2-dilation roughly 1
also.  The second retraction is accomplished by the map $id
\times \phi$.  The second retraction has 1-dilation roughly
$\epsilon^{-1}$ but it has 2-dilation 1.  To check the 2-dilation
of the retraction, we reason as follows.  The restriction of $id
\times \phi$ to $[0, \epsilon] \times A$ is the identity, and so
it has 2-dilation 1.  But the complement of $[0, \epsilon] \times
A$ in the domain of our map is just $\{ 0 \} \times A^c$.  Our
retraction maps this 2-dimensional set to the 1-dimensional set
$\{ 0 \} \times \partial A$.  Therefore, the second retraction
has 2-dilation 1. 

\vskip5pt
{\bf Lower bounds for the k-dilation}

Next we approach the problem from the other side, proving lower
bounds for the k-dilation $D_k(R,S)$.  Our lower bounds are based
on k-width and on the width-volume inequality.  Our estimates for
$D_k(R,S)$ depend on the dimensions of $R$ and $S$.  We adopt the
convention that $R$ and $S$ are n-dimensional rectangles, that $R$ has
dimensions $R_1 \le ... \le R_n$ and that $S$ has dimensions $S_1
\le ... \le S_n$.

The first lower bound on $D_k(R,S)$ comes from knowing the
k-width of $S$.  Suppose that $f$ is a degree 1 map from $R$ to
$S$ with k-diliation $\lambda$.  The rectangle $R$ can be sliced
into k-dimensional rectangles with dimensions $R_1 \times ...
\times R_k$, and these rectangles form a family of cycles
sweeping out $R$.  The image of each k-dimensional rectangle has
volume at most $\lambda R_1 ... R_k$.  The situation is
illustrated in Figure 3.

\includegraphics{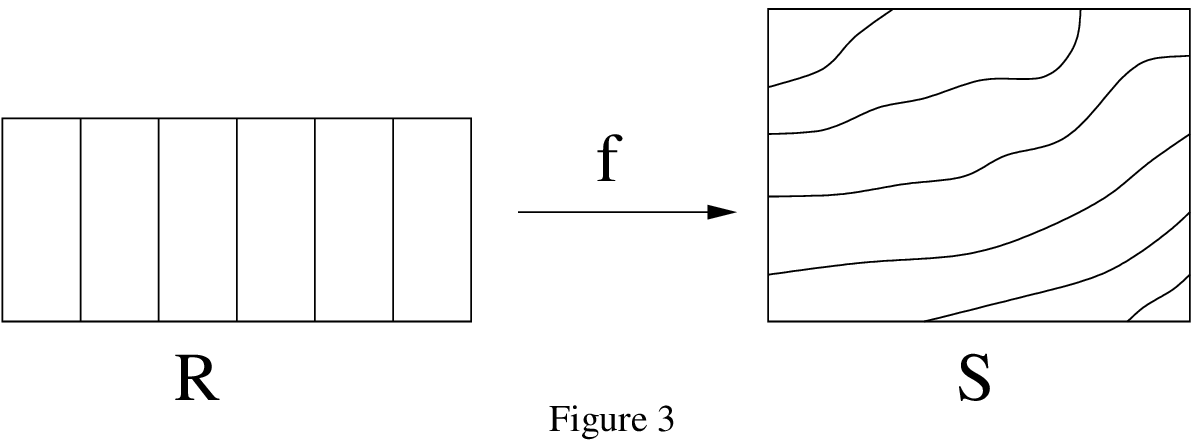}

The image of our family of rectangles is a family of k-cycles
sweeping out $S$.  According to Proposition 1, this family must
contain a cycle with volume at least $c(n) S_1 ... S_k$.  Since
each cycle in the family has volume at most $\lambda R_1 ...
R_k$, we get a lower bound for the k-dilation $\lambda$.

\begin{introprop} $D_k(R,S) > c(n) [S_1 ... S_k] / [R_1 ...
R_k]$.
\end{introprop}

We can get more estimates if, instead of considering the k-width
of $S$, we work with the k-width of subsets of $S$.  Let's see how
this idea works out in a particular case.  Suppose that $S$ is a
3-dimensional rectangle with dimensions $S_1 \le S_2 \le
S_3$.  Then $S$ contains many subrectangles with dimensions $S_1
\times S_2 \times S_2$.  We can find $N$ disjoint rectangles in $S$
with those dimensions, where $N$ is roughly $S_3/S_2$.  Call the
rectangles $V_i$.  Each one of these rectangles has 2-width
roughly $S_1 S_2$.  Now suppose that $f$ is a degree 1 map from
$R$ to $S$ with 2-dilation $\lambda$.  Then each of our
rectangles has a preimage $U_i = f^{-1}(V_i)$, and each of these
preimages has 2-width at least $S_1 S_2 / \lambda$.  The
situation is illustrated in Figure 4.

\vskip5pt

\includegraphics{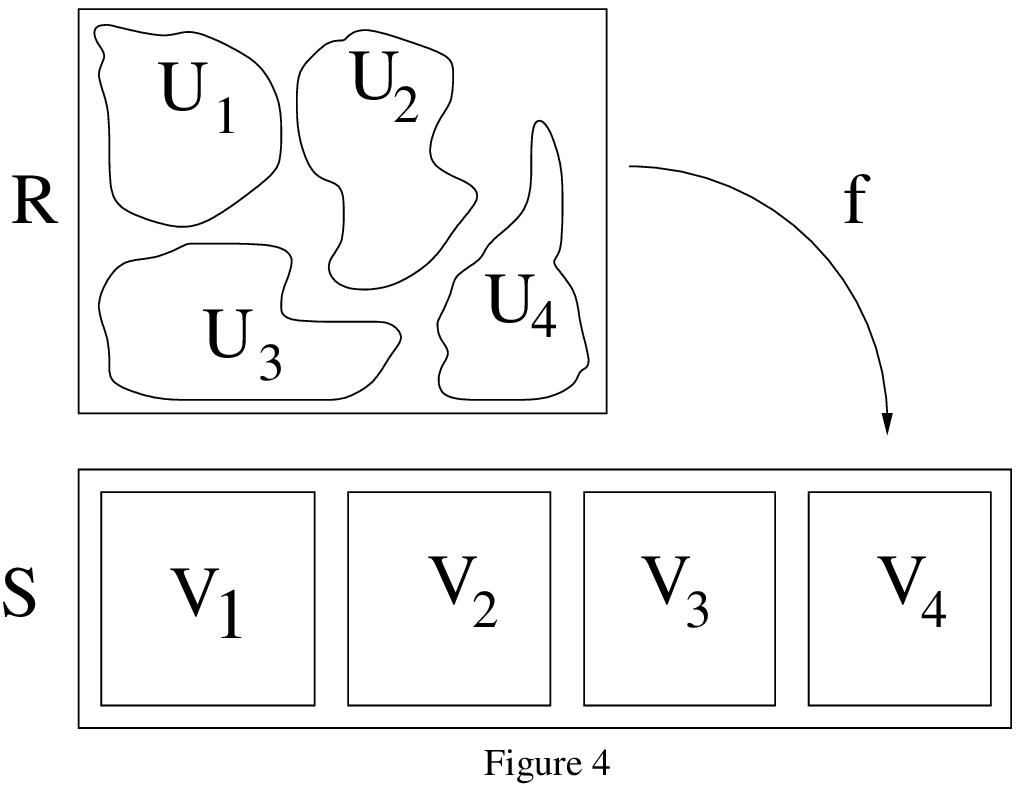}

We want to use this information to get a lower bound on
$\lambda$.  Since the sets $U_i$ are disjoint, one of them must
have volume at most $R_1 R_2 R_3 / N$.  We are led to the
following question: if $U \subset R$ is an open set with volume
$V(U)$, what is the largest possible 2-width of $U$?  Since $U$
is a subset of $R$, its 2-width is at most $R_1 R_2$.  For large
volumes $V(U)$, this upper bound is the best possible, but for
smaller volumes it can be improved.  Using the width-volume
inequality, we can bound the 2-width of $U$ by $C V(U)^{2/3}$. 
This upper bound is sharp for small volumes $V(U)$.  These upper
bounds can be improved if $V(U)$ is in the intermediate range
$R_1^3 << V(U) << R_1 R_2^2$.  An example of a set $U$ with
2-width roughly $V(U)^{2/3}$ is the round ball of volume $V(U)$,
which has radius roughly $V(U)^{1/3}$.  If $R_1^3 << V(U)$, then
this round ball does not fit in the rectangle $R$.  It turns out
that all subsets of $R$ with volume $V(U)$ are substantially
thinner than the round ball.  We make this precise in the
following proposition.

\begin{introprop} Let $R$ be a 3-dimensional rectangle with
dimensions $R_1 \le R_2 \le R_3$.  Suppose that $U \subset R$ is
an open set with volume $V(U)$.  Then the 2-width of U is at most
$C R_1^{1/2} V(U)^{1/2}$.
\end{introprop}

This estimate is a variation on the width-volume inequality
adapted to subsets of the rectangle $R$.  It improves on the
original inequality exactly when $R_1^3 << V(U)$.  The proof is
only a small modification of the proof of the width-volume
inequality.  Using this inequality to upper bound the 2-width of
one of the sets $U_i$, we get a new lower bound for $D_2(R,S)$.

\begin{introprop} If $R$ and $S$ are 3-dimensional rectangles
with dimensions $R_1 \le R_2 \le R_3$ and $S_1 \le S_2 \le S_3$,
then $D_2(R,S) > c [S_1 S_2^{1/2} S_3^{1/2}] / [R_1 R_2^{1/2}
R_3^{1/2}]$.
\end{introprop}

In the paper we carry out this idea for all values of k and n,
proving lower bounds for $D_k(R,S)$.  In the special case that
$k=n-1$, our lower bounds and the maps we will construct match up
well enough to determine $D_{n-1}(R,S)$ up to a constant factor.

\begin{introtheorem} Let $R$ and $S$ be n-dimensional
rectangles.  Let $R$ have dimensions $R_1 \le ... \le R_n$, and $S$
have dimensions $S_1 \le ... \le S_n$.  Let $Q_i$ denote the
quotient $S_i/R_i$.  Up to a constant factor $C(n)$, the optimal
$(n-1)$-dilation $D_{n-1}(R,S)$ is equal to the maximum of the following
list of n monomials in the variables $Q_i$.

The first n-1 monomials are given by $Q_1 ... Q_l (Q_{l+1} ...
Q_n)^{\frac{n-l-1}{n-l}}$, where l is an integer in the range $1
\le l \le n-1$.  The final monomial is $Q_2 ... Q_n$.
\end{introtheorem}

The algebra here is somewhat complicated, but the complicated
expressions in $Q_i$ are not the important point.  We have seen
that the snake map can have (n-1)-dilation much smaller than that
of the linear map.  For any two rectangles $R$ and $S$ we will
construct an explicit map with nearly optimal (n-1)-dilation. 
Depending on the rectangles, it may be a linear map, or it may be
a minor generalization of the snake map.  Up to a constant
factor, the expression in the theorem will turn out to be the
(n-1)-dilation of this map.  The lower bounds in the theorem
guarantee that the (n-1)-dilation of this map cannot be
substantially improved.  (On the other hand, for $2 \le k < n-1$,
the k-dilation of the snake map can be improved in some cases. 
For more information on this problem, see \cite{Gu1}.)

\vskip5pt
{\bf Related results and open questions}

The literature contains a couple of theorems in a similar spirit
to the width-volume inequality.  For example, in appendix 1 of
\cite{G1}, page 128, Gromov proved the following estimate connecting the
Uryson width and the area of a Riemannian 2-sphere.  (The Uryson
width is a different notion of width from the one in this paper. 
For a definition, see Gromov's book \cite{G2}, page 108.)

\begin{reftheorem} (Gromov) Let $(S^2, g)$ be a Riemannian
2-sphere with Uryson 1-width $W$ and area $A$.  Then $W < 2
A^{1/2}$.
\end{reftheorem}

\noindent Another geometric quantity related to the k-width is the
volume of the smallest stationary k-cycle in a Riemannian
manifold.  According to the work of Almgren \cite{A}, a closed
oriented Riemannian manifold $(M,g)$ contains a stationary
k-dimensional varifold with volume at most $W_k(M,g)$.  Recently,
Nabutovsky and Rotman proved several estimates for the length of
the shortest stationary 1-cycle in a Riemannian manifold.  One
important estimate is the following theorem from \cite{NR}.

\begin{reftheorem} (Nabutovsky, Rotman) A closed
Riemannian n-manifold $(M,g)$ of volume $V$ contains a stationary
1-cycle of length at most $C(n) V^{1/n}$.
\end{reftheorem}

\noindent Although these theorems are in a similar spirit to
Theorem 1, they don't give any upper bounds for k-width for any
value of k.  These theorems hold in a more general setting than
Theorem 1 because they apply to arbitrary Riemannian metrics,
whereas Theorem 1 applies only to domains in Euclidean space.

Comparing our result to the results of Gromov, Rotman, and
Nabutovsky, it seems reasonable to ask whether there is a
width-volume inequality for all Riemannian manifolds.  We phrase
this as a problem.

\begin{openprob} For which integers $k < n$ is there a constant
$C(k,n)$ so that for every closed oriented Riemannian n-manifold
$(M,g)$, the k-width is bounded in terms of the volume by the
formula $W_k(M,g) < C(k,n) \textrm{ Volume}(M,g)^{k/n}$?
\end{openprob}

\noindent In an appendix to this paper, we show that the answer
to this question is negative when $k = n-1$.  In other words, a
closed oriented Riemannian n-manifold may have volume 1 and
arbitrarily large (n-1)-width.  For $k < n-1$, the problem is
open.

Now we turn to some other open problems related to k-width.

\begin{openprob} If k is in the range $2 \le k \le n/2$, is it
possible to bound the linear k-width of a bounded open set $U
\subset \mathbb{R}^n$ in terms of its volume?
\end{openprob}

This problem is related to the problem of Besicovitch $(n,k)$
sets.  An $(n,k)$ set is defined to be a subset of $\mathbb{R}^n$
of Lebesgue measure zero containing a translate of every k-plane. 
Besicovitch gave examples of $(2,1)$ sets, and an easy
generalization shows that $(n,1)$ sets exist for all n.  On the
other hand, Falconer's theorem mentioned above proves that there
are no $(n,k)$ sets for $k > n/2$.  There is recent progress on
this problem starting with Bourgain's important paper \cite{Bo}
on the Kakeya maximal function.  Bourgain
proves that there are no $(4,2)$ sets and no $(7,3)$ sets, among
other results.  The problem of $(n,k)$ sets, however, is not
equivalent to the problem above.  For example, to prove that
there are no $(4,2)$ sets, Bourgain establishes the following
slightly weaker version of a linear width-volume inequality. 
(The theorem below follows immediately from Proposition 3.3 in
\cite{Bo}.)

\begin{reftheorem} (Bourgain) For each $\epsilon > 0$, there is a
constant $C_\epsilon$ so that the following estimate holds. 
Let $U$ be a bounded open subset in $\mathbb{R}^4$, with volume
$V(U)$ and diameter $D(U)$.  Then the linear 2-width of $U$ is
bounded by $C_\epsilon V(U)^{\frac{2-\epsilon}{4}} D(U)^\epsilon$.
\end{reftheorem}

\begin{openprob} Find the k-width of the unit n-cube.
\end{openprob}

The exact k-width of the unit n-sphere was determined by Almgren,
and in \cite{G1}, Gromov claims that the k-width of the unit
n-ball is exactly the volume of the unit k-ball.  Because a cube
contains a ball, the results of Almgren and Gromov show that the
k-width of the unit n-cube is bounded below by a constant
independent of n.  The linear k-width of the unit n-cube is known
to be 1 by a result of Vaaler \cite{V}.  Even this result is
difficult.  Another interesting problem along these lines,
described by Gromov in \cite{G1}, is to estimate the k-width of
the unit ball in the finite-dimensional Banach space
$l^\infty(n)$.  Gromov showed that an estimate independent of n
would lead to results in intrinsic Riemannian geometry.

\begin{openprob} (The sponge problem) Recall that an embedding
$I$ is called expanding if its derivative increases the length of
all vectors, or equivalently if it increases the lengths of all
curves.  For which dimensions n is there a constant $\epsilon >
0$ so that any bounded open set $U
\subset \mathbb{R}^n$ with volume less than $\epsilon$ admits an
expanding embedding into the unit n-ball?
\end{openprob}

A potential counterexample $U$ must have a small volume and a
large diameter.  It might resemble a sponge: a large cube from
which many tubes have been cut leaving a complicated region with
small volume.  The expanding embedding reminds me of squeezing
the water out of the sponge.  An affirmative answer to the sponge
problem would give a new proof of the width-volume inequality,
because the image of $U$ can be swept out by k-planes with volume
1, and the inverse images of these k-planes sweep out $U$ with
bounded volume.  I tried for a long time to construct the
expanding embeddings, but I wasn't able to do it even in
dimension 2.

\begin{openprob} Estimate $D_k(R,S)$ up to a constant factor, for
k in the range $2 \le k \le n-2$.
\end{openprob}

The cases $k=1$ and $k=n$ are elementary, and the case $k=n-1$ is
done in this paper.  This paper contains some lower bounds for
$D_k(R,S)$.  It also includes some degree 1 mappings which give
upper bounds for $D_k(R,S)$.  There are several more upper and
lower bounds for $D_k(R,S)$ in \cite{Gu1}.  The gap between the
best upper and lower bounds can be arbitrarily large, however. 
The work in \cite{Gu1} suggests that the cases $2 \le k \le n-2$
are a lot harder than the case $k = n-1$.  This specific problem
may not have any applications, but I think it's a good testing
ground to see how well we understand k-dilation.

\vskip5pt

{\bf The plan of the paper}

In the first section of the paper, we give the precise definition
of k-width and its basic properties.  We give Gromov's proof that
the k-width is not zero and estimate the k-widths of cubes and
rectangles.  In the second section of the paper, we carry out the
``bending planes around a skeleton'' construction and use it to
prove the width-volume inequality.  The second section also
includes the variation of the width volume inequality adapted to
subsets of a rectangle.  In the third section, we use these
estimates to prove lower bounds for the k-dilation of degree 1
maps.  In particular, we prove all the lower bounds in Theorem 2. 
In the fourth section, we construct the snake map and its
higher-dimensional analogues.  Using these maps, we check that
non-linear maps can have much smaller k-dilation than linear
maps.  Then we prove all the upper bounds in Theorem 2.  The
fifth section is an appendix which gives lower bounds for the
width of Riemannian manifolds using isoperimetric inequalities. 
It shows that a Riemannian manifold of volume 1 may have
arbitrarily large (n-1)-width.  The sixth section is a second
appendix, briefly explaining Falconer's bound on the linear
k-width for $k > n/2$.

Throughout the paper we use $c(n)$ and $C(n)$ to denote constants
that depend only on the dimension n.  The value of these
constants may change from line to line.  We use $C(n)$ to denote
a large constant and $c(n)$ to denote a small positive constant. 
As described above, when we talk about a rectangle $R$, we always
order its dimensions so that $R_1 \le ... \le R_n$.

This paper is based on a section of my thesis \cite{Gu1}, and I
would like to thank my advisor Tom Mrowka for his help and
support.  I am also grateful to the referees for their
constructive comments.

\section{The definition of k-width}

In this section, we will make precise the intuitive idea of
k-width described in the introduction.  Our first task is to give
a precise meaning to a family of k-cycles.  We will use families
of flat integral cycles which are continuous in the flat
topology.  Roughly speaking, the flat topology means that two
k-cycles are close to one another if their difference bounds a
(k+1)-chain with small volume.  The precise definition that we
give is somewhat technical.  A reader interested in the main
ideas of the paper might skip the definition and proceed with
only an intuitive idea of a family of cycles sweeping out M. 
After defining the k-width, we prove that it behaves
monotonically with respect to appropriate mappings and that it is
not zero.  At the end of the section, we estimate the k-widths of
cubes and rectangles.

Our explanation of the flat topology and flat cycles essentially
follows Fleming's paper \cite{F}.  We recall that an integral
Lipschitz k-chain in $(M,g)$ is a finite sum $\sum c_i f_i$, where
each $c_i$ is an integer and each $f_i$ is a Lipschitz map from
the k-simplex to M.  An integral k-chain is a special case of a
singular k-chain with integer coefficients.  We define the
boundary of a k-chain as in the singular homology theory.

We define the volume of the map $f_i$ to be the volume of the
k-simplex in the induced metric $f_i^*(g)$.  Then we define the
mass of the chain $\sum c_i f_i$ to be $\sum |c_i| \textrm{
volume}(f_i)$.  We abbreviate the mass of a k-chain C by $|C|$. 
The most important fact about k-chains and mass is the
isoperimetric inequality.  We will use the following rather weak
formulation of the isoperimetric inequality.

\begin{reftheorem} (Isoperimetric Inequality) For each compact
manifold $(M, g)$, there is a constant $\epsilon$ so that every
integral Lipschitz k-cycle with mass less than $\epsilon$ is
homologically trivial.  Moreover, if C is a k-cycle with mass $m
< \epsilon$, then $C = \partial D$ for some integral Lipschitz
$(k+1)$-chain D with mass less than $C(n) m$.
\end{reftheorem}

Now we define the flat norm.  The flat norm was introduced by
Whitney in \cite{W} and used to define spaces of cycles by
Fleming in \cite{F}.  The flat norm of a k-chain C is the infimal
value of $|C - \partial D| + |D|$ as D varies over all the
integral Lipschitz (k+1)-chains in M.  We define the flat distance
between chains $C_1$ and $C_2$ as the flat norm of $C_1 - C_2$. 
It may happen that the distance between two flat chains is zero. 
This occurs when the two chains parametrize the same geometric
object in different ways.  According to a result of Fleming (see
\cite{F}), the distance is zero only for this reason.  We do not
need this result however.  We identify any two chains separated
by flat distance zero.  The space of equivalence classes is now a
metric space.  The completion of this metric space is called the
space of integral flat k-chains in M.

We define the mass of a flat chain $C$ to be the infimal number $m$
so that there exists a sequence of integral Lipschitz chains $C_i$
converging to $C$ with mass less than $m$.

The boundary map on integral Lipschitz chains is bounded with
respect to the flat norm.  In fact, if C is a k-chain and D is a
(k+1)-chain, then the flat norm of $\partial C$ is at most $|C -
\partial D|$.  Taking the infimum over all D shows that the flat
norm of C is greater than or equal to the flat norm of the
boundary of C.  Therefore, we can define boundaries of flat
chains.  The flat k-cycles are the subset of flat k-chains with
boundary zero.  Let Z denote the space of flat k-cycles with the
flat topology.

By a family F of k-cycles in M, we mean a continuous map F from a
parameter space P to Z.  We will always assume that P is a finite
simplicial complex.  We define the width of the family F to be
$sup_{p \in P} |F(p)|$.

Our next task is to define what it means for a family of k-cycles
to sweep out M.  Morally, an i-dimensional family of k-cycles can
be glued together to form a (k+i)-cycle, but this is not
literally true for the space of flat cycles.  We now give a
construction that takes an i-dimensional family of cycles and
gives a (k+i)-cycle that, in some sense, is a small perturbation
of the family.

Let F be a family of k-cycles parametrized by P.  We take a
fine triangulation of P.  We pick a small number $\delta >
0$.  For each vertex v of the triangulation, we choose an integral
Lipschitz cycle $C(v)$, with flat distance less than $\delta$ from
$F(v)$, and mass less than $|F(v)| + \delta$.

(For completeness, we include the proof that such a cycle $C(v)$
exists.  By the definition of a flat cycle, we can take a
sequence of integral Lipschitz k-chains $C_i$ converging to $F(v)$
in the flat norm.  By the definition of mass, we may assume that
each chain $C_i$ has mass less than $|F(v)| + \epsilon$.  We have
to show that we can choose the chains $C_i$ to be cycles.  Because
the boundary operation is continuous in the flat norm, we know
that the flat norm of $\partial C_i$ converges to zero.  By the
definition of the flat norm, we may choose integral Lipschitz
k-chains $D_i$ so that $|\partial C_i - \partial D_i| + |D_i|
\rightarrow 0$.  According to the isoperimetric inequality, for
sufficiently large i, there is a k-chain $E_i$ with $\partial E_i
= \partial C_i - \partial D_i$ and $|E_i| < C |\partial C_i -
\partial D_i|$.  Now we define $\tilde C_i = C_i - D_i - E_i$. 
Each $\tilde C_i$ is an integral Lipschitz k-cycle.  Moreover,
$|D_i| + |E_i| \rightarrow 0$.  Since the mass controls the flat
norm, $\tilde C_i$ converges to $F(v)$ in the flat topology. 
Since $|D_i| + |E_i| \rightarrow 0$, the mass of $\tilde C_i$ is
less than $|F(v)| + 2 \epsilon$ for sufficiently large $i$.)

Now, since the triangulation is fine, we may assume that if
$v_1$ and $v_2$ are neighboring vertices, then the flat distance
between $C(v_1)$ and $C(v_2)$ is less than $3 \delta$.  By
definition, this means that there is an integral Lipschitz
$(k+1)$-chain D with $|C(v_1) - C(v_2) - \partial D| + |D|$ less
than $3 \delta$.  The k-cycle $C(v_1) - C(v_2) - \partial D$ must
have mass less than $3 \delta$.  According to the isoperimetric
inequality, it must bound a (k+1)-chain $D'$ with mass at
most $C(n) \delta$.  In other words, $C(v_1) - C(v_2) = \partial (D
+ D')$.  We know that the mass of $D + D'$ is bounded by $C(n)
\delta$.

Let E denote the edge from $v_1$ to $v_2$, oriented so that
$\partial E = v_1 - v_2$.  We define $C(E) = D + D'$.  We repeat
this operation for every edge of the triangulation of P.  For
each edge E with boundary $v_1 - v_2$, $C(E)$ is a (k+1)-chain
with boundary $C(v_1) - C(v_2)$ and mass at most $C(n) \delta$.

We continue this procedure inductively.  For each oriented
i-dimensional simplex $\Delta^i$ in P, we define a (k+i)-chain
$C(\Delta^i)$ with the following properties.  If the
boundary of the simplex $\Delta^i$ is equal to $\sum_j
\Delta_j^{i-1}$ with orientations, then the boundary of
$C(\Delta^i)$ is equal to $\sum_j C(\Delta_j^{i-1})$ as integral
Lipschitz cycles.  Moreover, $C(\Delta^i)$ has mass less than $C(n)
\delta$.  We can always choose such chains by using the
isoperimetric inequality in M, provided the initial number
$\delta$ is sufficiently small.

The map C taking oriented simplices of P to integral Lipschitz
chains in M can be thought of as a chain map between two chain
complexes.  The first complex is associated to P.  To describe
it, it is convenient to first pick an orientation for every
simplex in P.  The complex has i-chains consisting of sums $c_i
\Delta_i$, where $c_i$ is an integer and $\Delta_i$ is any
i-dimensional simplex of P.  The homology of this chain complex
is the simplicial homology of P with coefficients in
$\mathbb{Z}$.  The second chain complex consists of the integral
Lipschitz chains in M.  The homology of this complex is equal to the
singular homology of M with coefficients in $\mathbb{Z}$.  The
map C is a map from the first complex to the second complex,
taking $i$-chains to $(k+i)$-chains, and commuting with boundary
operations - in other words a chain map with shift k.  We call
such a chain map C a complex of k-cycles in M.

The chain map C induces a map from the simplicial homology of P
to the singular homology of M.  Since simplicial homology and
singular homology agree, we get a map from $H_i(P, \mathbb{Z})$
to $H_{k+i}(M, \mathbb{Z})$.  We call this map the gluing
homomorphism G.

The chain map C was not canonical.  On the contrary it involved
many choices.  Nevertheless, the gluing homomorphism does not
depend on these choices, as long as $\delta$ is sufficiently small.  To
see this, let $C_0$ and $C_1$ be two possible choices of chain
map following the construction above. Divide $P \times [0,1]$
into cells given by $\Delta \times \{0\},
\Delta \times \{1 \},$ and $\Delta \times [0,1]$, where $\Delta$
varies over the triangulation of P.  Now define $C(\Delta \times
\{0\}) = C_0(\Delta)$ and $C(\Delta \times \{1\}) = C_1(\Delta)$. 
Suppose that $\Delta^i$ is an i-simplex in P, and that the
boundary of $\Delta^i \times [0,1]$ is equal to $\Delta^i \times
\{ 1 \} - \Delta^i \times \{ 0 \} + 
\sum_j \Delta_j^{i-1} \times [0,1]$ with orientations. 
Proceeding inductively, we define $C(\Delta^i \times [0,1])$ to
be a (k+i+1)-chain with boundary $C(\Delta^i \times \{1\}) -
C(\Delta^i \times \{0\}) + \sum_j C(\Delta_j^{i-1} \times
[0,1])$.  Again, as long as we assume $\delta$ sufficiently
small, we can construct all these chains using the isoperimetric
inequality.  We can view C as a chain map as well, and it induces
a map on homology from $H_i(P \times [0,1],
\mathbb{Z})$ to $H_{k+i}(M, \mathbb{Z})$.  This map agrees with
the gluing homomorphisms induced by both $C_0$ and $C_1$, and so
these two gluing homomorphisms agree with one another.

We say that a family of k-cycles F sweeps out M if its gluing
homomorphism $G: H_{n-k}(P, \mathbb{Z}) \rightarrow H_n(M,
\mathbb{Z})$ is non-trivial.

We define the k-width of M to be the infimal width of any
family of k-cycles sweeping out M.  We denote the k-width of M by
$W_k(M)$.

(On a manifold with boundary, the k-width is defined using
relative flat k-cycles.  The above arguments also apply 
to relative chains and cycles.  In this case, the gluing
homomorphism maps $H_i(P, \mathbb{Z})$ to $H_{k+i}(M, \partial M,
\mathbb{Z})$.  We say that a family F of relative k-cycles sweeps
out M if the gluing homomorphism from $H_{n-k}(P)$ to $H_n(M,
\partial M)$ is non-trivial.)

Next we discuss some basic properties of the k-width.  

\vskip10pt

\centerline{\bf Monotonicity}

The k-width is monotonic in two respects.  First, if M is an open
subset of $(N,g)$ with the induced metric, then the k-width of M
is at most the k-width of N.  Let F be a family of k-cycles
sweeping out N with width less than $W_k(N) + \epsilon$.  There
is a restriction map which takes integral Lipschitz chains in N
to integral Lipschitz chains in M.  This map commutes with the
boundary action and is bounded in the flat norm.  Therefore, it
takes flat k-cycles in N to (relative) flat k-cycles in M. 
Restricting F to M gives a family of cycles sweeping out M with
width less than $W_k(N) + \epsilon$.  This proves the first form
of monotonicity.

The second form of monotonicity concerns maps from M to N.  If f
is a Lipschitz map between compact Riemannian manifolds M and N,
then it maps integral Lipschitz chains in M to integral Lipschitz
chains in N.  This induced map is continous with respect to the
flat distance, so it maps flat k-chains in M to flat k-chains in
N.  If M and N are both closed, it maps flat k-cycles in M to
flat k-cycles in N.  If M and N both have boundaries, and if f
maps $(M, \partial M)$ to $(N, \partial N)$, then f maps relative
flat k-cycles in M to relative flat k-cycles in N.  If M has
boundary and N is closed, and if f maps $(M, \partial M)$ to $(N,
*)$, where $*$ is a point of N, then f takes relative flat
k-cycles in M to flat k-cycles in N.  The construction above
shows that the gluing homomorphism is natural.

If f is a Lipschitz map from M to N with Lipschitz constant L and
non-zero degree then $L^k W_k(M) \ge W_k(N)$.  To see this, let F
be a family of flat k-cycles sweeping out M with width less than
$W_k(M) + \epsilon$.  Since f has non-zero degree, the image
$f(F)$ is a family of k-cycles sweeping out N.  Since f maps each
integral Lipschitz k-chain with mass M to one with mass at most
$L^k M$, the width of $f(F)$ is less than $L^k(W_k(M) +
\epsilon)$.

The last estimate really only depended on how much the map f
stretched k-dimensional volumes.  Let f be a piecewise smooth
map.  Recall that the k-dilation of f is at most $\lambda$ if f
maps each k-dimensional submanifold of the domain with volume V
to an image with volume at most $\lambda V$.  For more
information on k-dilation, see the beginning of section 3.  Let f
be a piecewise smooth map of non-zero degree from M to N with
k-dilation $\lambda$.  Since Lipschitz maps can be
well-approximated by $C^1$ maps, f takes each integral Lipschitz
chain with mass M to an integral Lipschitz chain with mass at most
$\lambda M$.  Let F be a family of k-cycles sweeping out M with
width at most $W_k(M) + \epsilon$.  Then f(F) sweeps out N with
width at most $\lambda (W_k(M) + \epsilon)$.  Therefore, $\lambda
W_k(M) \ge W_k(N)$.

\vskip10pt

\centerline{\bf Non-degeneracy}

The non-degeneracy property says that $W_k(M) > 0$ for any
$(M,g)$.  Because of the monotonicity estimates for $W_k$, the
non-degeneracy follows for every n-manifold as soon as we know
that $W_k(S^n) > 0$ for the standard round metric on $S^n$. 
Gromov gave an elegant elementary proof of this fact which we
include here in our language.  This proof originally appeared in
\cite{G1}.

\begin{prop} (Gromov) For the standard round metric on $S^n$, the
k-width $W_k(S^n)$ is greater than $c(n) > 0$.
\end{prop}

\proof Suppose not.  Then there is a family F of k-cycles
sweeping out $S^n$ with width less than $\epsilon$.  Above, we
constructed a complex of cycles C based on F, with a non-trivial
gluing map.  For each vertex v of the triangulation of the
parameter space P, we had $|C(v)| < \epsilon + \delta$, and for
each higher dimensional simplex $\Delta^i$ of P, we had
$|C(\Delta^i)| < C(n) \delta$.  These estimates hold for a number
$\delta$ as small as we like.  We assume $\delta$ much smaller
than $\epsilon$.  

Using the isoperimetric inequality, each cycle $C(v)$ can be
filled by a (k+1)-chain of mass less than $C(n)
\epsilon$.  For each vertex v of the triangulation, define $Fill(v)$
to be such a filling.  Now let E be an oriented edge of the
triangulation with $\partial E = v_1 - v_2$.  We define $\bar
C(E) = C(E) - Fill(v_1) + Fill(v_2)$.  Since the boundary of $Fill(v_i) =
C(v_i)$, $\bar C(E)$ is a (k+1)-cycle.  The mass of $\bar
C(E)$ is bounded by $C(n) \epsilon$.  Next, using the isoperimetric
inequality again, choose an oriented (k+2)-chain $Fill(E)$ with
boundary $\bar C(E)$, and with mass bounded by $C(n)
\epsilon$.

We then repeat this construction for the higher-dimensional
simplices in the triangulation of P, working one skeleton at a
time.  For each i-simplex $\Delta^i$ of P, we define a
(k+i)-cycle $\bar C(\Delta^i)$ and a $(k+i+1)$-chain
$Fill(\Delta^i)$.  They have the following properties.

1. If $\Delta^i$ is an i-simplex of P, and the boundary of $\Delta^i
= \sum_j \Delta_j^{i-1}$, then

$\bar C(\Delta^i) = C(\Delta^i) - \sum_j Fill(\Delta_j^{i-1}).$

2. If $\sum_l c_l \Delta^i_l$ is any i-cycle in P, then

$\sum_l c_l C(\Delta^i_l) = \sum c_l \bar C(\Delta^i_l).$

3. The boundary of $Fill(\Delta^i) = \bar C(\Delta^i)$.

4. The mass of each $C(\Delta^i)$ and each $Fill(\Delta^i)$ is
bounded by $C(n) \epsilon$.

To prove that we can find $\bar C(\Delta^i)$ and $Fill(\Delta^i)$
we work inductively.  We already did the case $i=1$, which
anchors the induction.  We assume that the above properties hold
for simplices of dimension at most i-1.  We define $\bar
C(\Delta^i)$ by using the formula in 1.  We have to check that
$\bar C(\Delta^i)$ is a cycle.  Its boundary is $\sum_j
C(\Delta_j^{i-1}) - \sum_j \bar C(\Delta_j^{i-1})$.  According to
the equation in 2 and the inductive hypothesis, this expression
vanishes.  Next, we have to check that $\bar C$ obeys equation 2
for i-simplices.  Let $\sum_l c_l \Delta^i_l$ be an i-cycle in P. 
Let the boundary of $\Delta^i_l$ be $\sum_j \Delta^{i-1}_{l,j}$. 
Because $\sum_l c_l \Delta^i_l$ is an i-cycle, $\sum_{l,j} c_l
\Delta^{i-1}_{l,j} = 0$.  Now $\sum_l c_l \bar C(\Delta^i_l) =
\sum_l (c_l C(\Delta^i_l) + c_l \sum_j
Fill(\Delta^{i-1}_{l,j}))$, and the terms in the second sum
cancel because $\sum_l c_l \Delta^i_l$ is a cycle.  Therefore,
$\bar C$ obeys equation 2 for i-simplices.  Because $C(\Delta^i)$
has mass at most $\delta$ and $Fill(\Delta^{i-1}_j)$ has mass at
most $C(n) \epsilon$, $\bar C(\Delta^i)$ has mass at most $C(n)
\epsilon$.  Therefore, we can use the isoperimetric inequality to
define $Fill(\Delta^i)$ with mass at most $C(n) \epsilon$.

Since F sweeps out $S^n$, there must be an (n-k)-cycle $a$ in P
with $G(a)$ non-trivial in $H_n(S^n, \mathbb{Z})$.  Write $a = \sum c_m
\Delta_m^{n-k}$.  By definition $\sum c_m C(\Delta_m^{n-k})$ has
a non-trivial homology class in $H_n(S^n, \mathbb{Z})$.  But this
sum is equal to a sum of cycles $\sum c_m \bar C(\Delta_m^i)$. 
Each cycle $\bar C(\Delta_m^i)$ has mass less than $C(n) \epsilon$,
and hence is null-homologous in $S^n$.  This contradiction
finishes the proof. \endproof

We now estimate the k-widths of some simple shapes.  The k-width
of the unit n-cube is at most 1, because it is swept out by
parallel k-planes each meeting it in a unit k-cube.  Because of
the non-degeneracy proposition and the monotonicity estimate, the
k-width of the unit n-cube is at least $c(n) > 0$.  

Applying the monotonicity estimate to our bounds for the unit
cube, we can estimate the k-width of any rectangle.  Let $R$ be a
rectangle with dimensions $R_1 \le ... \le R_n$.  In other words,
R is the product $[0, R_1] \times ... \times [0,R_n]$.

\begin{prop} The k-width of the rectangle $R$ is roughly $R_1 ...
R_k$.  More precisely, $c(n) R_1 ... R_k \le
W_k(R) \le R_1 ... R_k$.
\end{prop}

\proof To get the upper bound, simply consider the projection of
R onto the last (n-k) coordinates.  Each fiber of this projection
is a k-dimensional rectangle with volume $(R_1 ... R_k)$, and the
fibers fit together to form a family of k-cycles sweeping out R.

To get the lower bound, consider the map from $R$ to the unit cube
which sends $(x_1, ..., x_n)$ to $(x_1 / R_1 ,... , x_n/R_n)$. 
This map has degree 1 and k-dilation $(R_1 ... R_k)^{-1}$. 
Therefore $(R_1 ... R_k)^{-1} W_k(R) \ge W_k(C)$, where C denotes
the unit cube.  According to Proposition 1.1, $W_k(C) \ge c(n) > 0$,
so $W_k(R) \ge c(n) R_1 ... R_k$. \endproof

\section{The width-volume inequality}

In this section, we prove Theorem 1.

\begin{theorem} (Width-volume inequality) Let $U$ be a bounded open
set in $\mathbb{R}^n$
with volume $V(U)$ and k-width $W_k(U)$.  Then $W_k(U) < C(n)
V(U)^{k/n}$.
\end{theorem}

\proof By a scaling argument, it suffices to prove the theorem when the
volume of $U$ is 1.

The first step in the proof is to translate the unit lattice so
that its k-skeleton meets $U$ in a region of controlled volume. 
Let $S(x)$ denote the k-skeleton of the unit cubical lattice
centered at x, with axes parallel to the coordinates.  Since the
volume of $U$ is 1, the average volume of $U \cap S(x)$ as x varies
over the unit cube is equal to $n \choose k$.  We can choose a
point x so that the volume of $U \cap S(x)$ is no more than the
average value $n \choose k$.  From now on, we refer to S(x)
simply as S.  

The second step in the proof is to construct a family of cycles
sweeping out U, each of which lies mostly in the k-skeleton S.

\begin{construction} (Bending planes around a skeleton) Let
$B(R)$ denote the ball of radius $R$ in $\mathbb{R}^n$, and let S
be the k-skeleton of a unit lattice.  Then there is a family F of
k-cycles sweeping out $B(R)$ with the following properties.  Each
cycle in F lies in $S$ except for a subset of mass less than $C(n)$. 
A cycle in F may contain some portions of $S$ with multiplicity
greater than 1, but this multiplicity is bounded by $C(n)$.  (The
constant C(n) depends only on n; it does not depend on the radius
R.)
\end{construction}

Using this family of k-cycles we finish the proof of Theorem 1. 
By choosing $R$ sufficiently large, we may assume that $U$ lies in
the ball $B(R)$.  Then we consider the restriction of F to U,
which is a family of k-cycles sweeping out U.  To prove Theorem
1, we only have to check that each k-cycle in this family has
mass at most $C(n)$.  Let $E$ be a k-cycle in F.  In other words,
we have to check that the mass of E restricted to $U$ is at most
$C(n)$.  We divide the restriction of $E$ to $U$ into two pieces. 
We let $E_1$ be the part of this restriction which is contained
in S, and we let $E_2$ be the part of this restriction which is
not contained in S.  The chain $E_1$ is contained in S, and
according to Construction 1, it has multiplicity at most $C(n)$. 
Therefore its mass is at most $C(n) |S \cap U| \le {n \choose k}
C(n)$.  On the other hand, according to Construction 1, $E_2$ has
mass at most $C(n)$.  \endproof

Next we turn to the proof of Construction 1.

\proof We begin with a family of parallel k-planes.  Let P be an
$(n-k)$-plane through the origin, in general position with
respect to S.  Let $F_0$ be the family of all k-planes
perpendicular to P.  To bend the planes, we will construct a
degree 1 proper PL map $\Psi$ from $\mathbb{R}^n$ to itself.  Our
family F will be $\Psi(F_0)$.  In other words, the cycles in F
will be $\Psi(Q)$ as Q varies over all the k-planes perpendicular
to P.  Because $\Psi$ is degree 1, this family of cycles sweeps
out the ball $B(R)$ or any other open set.

The reader can roughly imagine $\Psi$ as follows.  Let T denote
the dual $(n-k-1)$-skeleton to S, and let $T_\epsilon$ denote the
$\epsilon$-neighborhood of T.  The mapping $\Psi$ retracts the
complement of $T_\epsilon$ onto $S$ while stretching $T_\epsilon$
to fill all of $\mathbb{R}^n - S$.

The idea of our proof is as follows.  Since $\Psi$ retracts the
complement of $T_\epsilon$ onto S, the map $\Psi$ takes $Q \cap
(\mathbb{R}^n - T_\epsilon)$ into S.  On the other hand, we will
try to control the size of $Q \cap T_\epsilon$ using the fact
that T is $(n-k-1)$-dimensional and Q is $k$-dimensional. 
Because of these dimensions, a generic plane Q will not intersect
T at all.  The set of planes Q in $F_0$ which intersect T has
codimension 1.  The set of planes which intersect T twice has
codimension 2, and so on.  Therefore, each plane Q intersects T
at most $(n-k)$ times.  Using this kind of argument, we will show
that $Q \cap T_\epsilon$ is contained in a union of $(n-k)$ small
balls.  Finally we will have to analyze the action of $\Psi$ on
each of these small balls.  This last step requires us to write
down the map $\Psi$ carefully.

We state the properties of the map $\Psi$ we will use in the form
of a lemma.

\begin{lemma} For each $\epsilon > 0$, there is a
piecewise-linear map $\Psi$ from $\mathbb{R}^n$ to itself with
the following properties.  The map $\Psi$ is linear on each
simplex of a certain triangulation of $\mathbb{R}^n$.  Each
top-dimensional simplex of this triangulation is labelled good or
bad.  For each good simplex $\Delta$, $\Psi(\Delta)$ lies in S. 
Each bad simplex lies in $T_\epsilon$.  The triangulation and the
map obey the following bounds.

1. The number of simplices of our triangulation meeting any unit
ball is bounded by $C(n)$.

2. The displacement $|\Psi(x) - x|$ is bounded by $C(n)$.

3. The diameter of each simplex is bounded by $C(n)$.
\end{lemma}

The only tricky part in checking this lemma is to get the bounds
with constants that don't depend on $\epsilon$.  We defer the
proof of the lemma until we finish the construction.

Because of the displacement bound, the map $\Psi$ is proper.  We
can deform $\Psi$ to the identity by taking $\Psi_t (x) = (1-t)
\Psi(x) + t x$.  Then $\Psi_0$ is equal to $\Psi$ and $\Psi_1$ is
the identity.  Each map $\Psi_t$ also obeys the displacement
bound, so they are all proper.  Therefore, $\Psi$ has degree 1.

We think of the cycle $\Psi(Q)$ as a sum of chains $\sum \Psi(Q
\cap \Delta)$ as $\Delta$ varies over all the simplices of our
triangulation.

We first consider the contribution to $\Psi(Q)$ coming from the
good simplices.  For each good simplex $\Delta$, $\Psi(\Delta)$
lies in $S$, and so $\Psi(Q \cap \Delta)$ also lies in $S$. 
Since $\Psi$ is linear on $\Delta$, the image $\Psi(Q \cap
\Delta)$ has multiplicity at most 1.  Next, we bound the
multiplicity of the sum $\sum_\Delta \Psi(Q \cap \Delta)$ as
$\Delta$ ranges over all the good simplices.  Because of the
displacement bound in Lemma 2.1, the multiplicity of this sum at
a point $s$ in $S$ only depends on the contributions from good
simplices $\Delta$ in a ball around $s$ of radius $C(n)$.  But
estimate 1 in Lemma 2.1 tells us that there are less than $C(n)$
simplices in this ball.

We now consider the contribution to $\Psi(Q)$ coming from the bad
simplices.  Since we are only proving bounds for the restriction
of $\Psi(Q)$ to the ball $B(R)$, we only need to consider the bad
simplices $\Delta$ so that $\Psi(\Delta)$ intersects the ball
$B(R)$.  Because of the displacement bound in Lemma 2.1, we only
need to consider the bad simplices $\Delta$ in the ball of radius
$R + C(n)$.  We let $B(R')$ denote this larger ball, and from now
on we only consider the bad simplices in this ball.  This
argument about balls is not the main point, but it is technically
easier to proceed this way because we only have to consider
finitely many simplices. 

Next we show that if $\epsilon$ is sufficiently small, the plane
Q intersects less than $C(n)$ bad simplices.  This estimate is
the main idea of the proof.  Recall that Q is a plane
perpendicular to the (n-k)-plane P.  Let $\pi$ denote the
orthogonal projection from $\mathbb{R}^n$ onto P.  The plane Q is
one of the fibers of $\pi$.  Let $T'$ denote the finite complex
$T \cap B(R')$.  Note that $T'$ is contained in a finite union of
$(n-k-1)$-planes.  We denote these planes as $T_i$.  Since P is
in general position with respect to T, the projections $\pi(T_i)$
are a finite set of $(n-k-1)$-planes in P, meeting transversely. 
Therefore, any point $p$ in P lies in $\pi(T_i)$ for at most
$(n-k)$ values of $i$.  Since the number of planes is finite, we
can choose $\epsilon$ sufficiently small so that any point $p$ in
$P$ lies within $\epsilon$ of $\pi(T_i)$ for at most $(n-k)$ values
of $i$.  Since $Q$ is a fiber of $\pi$, it meets the
$\epsilon$-neighborhood of $T_i$ for at most $(n-k)$ values of
$i$.  Since $Q$ is transverse to each of these planes, $Q
\cap T'_\epsilon$ is contained in $(n-k)$ balls each of radius
$C(n) \epsilon$.  Because of estimate 1 in Lemma 2.1, these balls
meet at most $C(n)$ bad simplices, and so $Q$ intersects at most
$C(n)$ bad simplices.

Finally, we bound the volume of $\Psi(Q \cap \Delta)$ where
$\Delta$ is a bad simplex.  Because of the diameter bound and the
displacement bound, $\Psi(\Delta)$ has diameter at most $C(n)$. 
Therefore, $\Psi(Q \cap \Delta)$ is a portion of k-plane with
diameter at most $C(n)$.  We conclude that $\Psi(Q \cap \Delta)$
has volume at most $C(n)$.

We finish the proof of Construction 1 by assembling these
estimates.  The cycle $\Psi(Q)$ lies in $S$ except for the
contributions from the bad simplices.  There are at most $C(n)$
bad simplices, and each bad simplex contributes mass at most
$C(n)$, and so the cycle $\Psi(Q)$ lies in $S$ except for a
portion with mass at most $C(n)$.  The multiplicity of $\Psi(Q)$
is also bounded.  We already bounded the contribution to the
multiplicity coming from the good simplices.  Since Q intersects
only $C(n)$ bad simplices in $B(R')$, the contribution to the
multiplicity coming from the bad simplices is also bounded.  \endproof

Before we go on, let us clarify which constants depend on which
other constants.  The most important point is that the constants
$C(n)$ depend only on n.  The constants $C(n)$ don't depend on R. 
On the other hand, the size of $\epsilon$ that we need to make
the above construction work does depend on R.  Therefore, we need
to prove Lemma 2.1 with constants $C(n)$ that don't depend on
$\epsilon$.  We now give the proof of Lemma 2.1.

\proof We will now construct the map $\Psi$. 
We begin by constructing the triangulation of good and bad
simplices.  First we need to make some definitions.  If A is a
k-dimensional face in S, then we define the link of A in the
following way.  The set A is defined by equations $x_i = a_i$ for
(n-k) coordinates i, and equations $a_j \le x_j \le a_j + 1$ for
the other k coordinates.  There is an (n-k) cube transverse to A
given by the equations $a_i - 1/2 \le x_i \le a_i + 1/2$ for the
(n-k) coordinates i above, and $x_j = a_j + 1/2$ for the other k
coordinates.  This cube is simply the (n-k) cube centered at the
center of A, perpendicular to A, with axes parallel to the
coordinate axes.  The link of A is defined to be the boundary of
this (n-k)-cube.  It consists of $2(n-k)$ (n-k-1)-cubes, each of
which is an (n-k-1)-dimensional face of T.  If B is an (n-k-1)
dimensional face of T, we define the link of B in an analogous
way.  It is a topological k-sphere consisting of 2(k+1)
k-dimensional faces of S.  We let A denote a k-dimensional face
of $S$ and B an (n-k-1)-dimensional face of T.  A quick calculation
shows that A is in the link of B if and only if B is in the link
of A.  For each pair (A, B) of faces with A in the link of B, we
define $K(A,B)$ to be the convex hull of the union of A and B. 

Next we check that the sets $K(A,B)$ tile $\mathbb{R}^n$.  The
hyperfaces of the tile $K(A,B)$ correspond to pairs $(A, b)$
where b is an (n-k-2)-face in the boundary of B, or pairs $(a,
B)$, where a is a (k-1)-face in the boundary of A.  (The
corresponding face is just the convex hull of A and b, or of a
and B.)  Each face borders exactly two tiles in our tiling. 
Given a face $(A, b)$, let $B'$ be the (n-k-1)-face in the link
of A which lies on the other side of b from B.  Then $K(A, B')$
is the only other tile with $(A, b)$ as a face.  Therefore, the
tiles form a pseudo-manifold, and the embedding of the tiles is
an orientation preserving proper map from the tile space to
$\mathbb{R}^n$.  The intersection of $K(A,B)$ with the skeleton S
is equal to A.  In particular, the only tiles that come near to
the center of A are tiles $K(A,B)$ for some B in the link of A. 
It is easy to check that a typical point close to the center of A
lies in exactly one of the tiles $K(A,B)$.  Therefore, the tiles
have disjoint interiors and cover all of space.

Any two tiles in our tiling are isometric.  After renumbering the
coordinates, translating, and reflecting, we can assume that A
and B have the following simple form.  The face A is given by the
inequalities $0 \le x_i \le 1$ for i from 1 to k, $x_i = 0$ for i
from $k+1$ to $n$.  The face B is given by inequalities $-1/2 \le
x_i \le 1/2$ for i from $k+1$ to $n-1$, the equations $x_i = 1/2$
for i from 1 to k, and $x_n = 1/2$.  The convex set $K(A,B)$ is
given by the inequalities $0 \le x_n \le 1/2$, $-x_n \le x_i \le
x_n$ for i from $k+1$ to $n-1$, and $|1/2 - x_i| \le |1/2 - x_n|$
for i from $1$ to $k$.  

We now divide each tile $K(A,B)$ into good and bad parts.  The
good part of $K(A,B)$ is given by $x_n \le 1/2 -
\epsilon$ and denoted $K_G(A,B)$.  The bad part of $K(A,B)$ is
given by $x_n \ge 1/2 - \epsilon$ and denoted by $K_B(A,B)$.  In
other words, the bad part of $K(A,B)$ lies in a small
neighborhood of B, and its complement is the good part.  Since B
is a face of T, the bad part of $K(A,B)$ lies in $C(n) \epsilon$
neighborhood of T.

If K is any convex polyhedron, we can define a barycentric
triangulation for K as follows.  For each face F of K, of any
dimension, let $c(F)$ denote the center of mass of F.  The
triangulation of the 0-skeleton of K is trivial.  Now suppose we
have triangulated the i-skeleton of K.  We extend this
triangulation to each $(i+1)$-face F of K, by taking the cone
from $c(F)$ to the triangulation on the boundary of F.  A good
thing about the barycentric triangulation is that if two convex
polyhedra intersect in a face of any dimension, then the two
barycentric sub-divisions of that face agree.  Therefore,
applying the barycentric subdivision to each good and bad
polyhedron in our tiling, we get a triangulation of
$\mathbb{R}^n$.  This is the triangulation that appears in the
statement of the lemma.  The map $\Psi$ will be linear on each
simplex of this triangulation.  We call a top-dimensional simplex
good if it lies in the good part of $K(A,B)$ and bad if it lies
in the bad part of $K(A,B)$.

At this point, we can check some of the bounds in the lemma.  The
number of simplices in a unit ball is bounded by $C(n)$.  The
number of tiles $K(A,B)$ does not depend on $\epsilon$ at all. 
The combinatorial structure of the tiles $K_G(A,B)$ and
$K_B(A,B)$ also does not depend on $\epsilon$.  Therefore, the
number of simplices in the barycentric triangulation also does
not depend on $\epsilon$.  Each simplex in contained in some set
$K(A,B)$ and so has diameter at most $C(n)$.  Also, each bad
simplex lies within a $C(n) \epsilon$-neighborhood of T. 

To finish the proof of the Lemma, we need to construct the map
$\Psi$.  We will have to check that $\Psi$ is linear on each
simplex of our triangulation, that $\Psi$ maps each good simplex
into S, and that $\Psi$ obeys the displacement bound.

The map $\Psi$ will take $K_G(A,B)$ onto A and $K_B(A,B)$ onto
$K(A,B)$.  We will specify the value of $\Psi$ at the center of
each face of $K_G(A,B)$ and of $K_B(A,B)$.  We then define $\Psi$
to be the unique function which is linear on the barycentric
subdivision and takes these values at the centers of the faces. 
To carry this out, we must write down all of the faces in
$K_G(A,B)$ and $K_B(A,B)$.

The faces in $K(A,B)$ are as follows.  First, any face a of A. 
Second, any face b of B.  Third, the convex hull of any face a of A
and any face b of B.  (These faces may have any dimensions.)

The faces of $K_G(A,B)$ are as follows.  First, any face a of A. 
Second, the intersection of the convex hull of a and b with the
set $x_n \le 1/2 - \epsilon$.  Third, the intersection of the
convex hull of a and b with the set $x_n = 1/2 - \epsilon$.  In
each case, we define $\Psi(c(F))$ to be the center of a. 
Therefore, $\Psi$ maps $K_G(A,B)$ into A.

The faces of $K_B(A,B)$ are as follows.  First, any face b of B. 
Second, the intersection of the convex hull of a and b with the
set $x_n \ge 1/2 - \epsilon$.  Third, the intersection of the
convex hull of a and b with the set $x_n = 1/2 - \epsilon$.  In
the first case, we define $\Psi(c(F))$ to be the center of b.  In
the last two cases, we define $\Psi(c(F))$ to be the center of a.

If a certain face F belongs to several different polyhedra, then
we have to check that our definition for $\Psi(c(F))$ is
consistent.  For this purpose, it suffices to check that the
faces a and b are defined consistently.  The face b is recovered
as the largest face of T in the $n \epsilon$ neighborhood of F,
provided $\epsilon$ is sufficiently small.  The face a is
recovered as the smallest face of $S$ so that F is contained in the
convex hull of a and b.

From the construction, we see that $\Psi$ is linear on each
simplex of our triangulation.  If $\Delta$ is a good simplex,
then each vertex of $\Delta$ corresponds to the center of a face
of $K_G(A,B)$, and so it gets mapped to a point in A.  Since A is
convex, the simplex $\Delta$ is mapped into A, and so
$\Psi(\Delta)$ lies in S.  Finally, if $\Delta$ denotes any good
or bad simplex in $K(A,B)$, then $\Psi(\Delta)$ lies in $K(A,B)$,
and so $\Psi$ obeys the displacement bound. \endproof

There is an analogue of Theorem 1 for the widths of functions
instead of sets.  Once we define the k-width of a function, the
proof is exactly the same.  Let $f$ be a compactly supported
function on $\mathbb{R}^n$ which is greater than or equal to
zero.  Let $B(R)$ denote a large ball containing the support of
f.  If F is a family of k-cycles in $B(R)$, then we define the
k-width of F to be the supremum of $\int_{F(p)} f$ over all p in
the parameter space of F.  We define the k-width of f to be the
infimal W so that there is a family of k-cycles sweeping out $B(R)$
with k-width less than W.  The k-width of f is denoted $W_k(f)$.

\begin{prop} If f is a function with compact support on
$\mathbb{R}^n$, and $0 \le f \le 1$, then
$W_k(f) < C(n) (\int f)^{k/n}$.
\end{prop}

\proof After rescaling the coordinates, it suffices to prove that
$W_k(f) < C(n)$ when $\int f = 1$.

By translating the k-skeleton of the unit lattice, we can arrange
that $\int_S f \le {n \choose k}$.  Next we apply Construction 1,
bending planes around the skeleton S.  This construction gives us
a family F of k-cycles sweeping out $B(R)$.  To prove the
proposition, we have to bound the integral $\int_{E} f$, where
$E$ is a k-cycle in the family F.  We define $E_1$ to be the part
of $E$ which is contained in $S$, and we define $E_2$ to be the
part of $E$ which is not contained in S.  The chain $E_1$ is
contained in $S$ and, according to Construction 1, it has
multiplicity at most $C(n)$.  Therefore, $\int_{E_1} f \le C(n)
\int_S f \le C(n) {n \choose k}$.  On the other hand, $E_2$ has
mass at most $C(n)$.  Because $0 \le f \le 1$, we have the
bound $\int_{E_2} f \le C(n)$. \endproof

In our applications, we will need a width-volume inequality
adapted to subsets of rectangles.  Let $R$ be the n-dimensional
rectangle with dimensions $R_1 \le ... \le R_n$.  Let $U$ be a
subset of $R$ with volume $V(U)$.  What is the largest possible
k-width of such a set $U$?  We already know that $W_k(U) < C
V(U)^{k/n}$, but this estimate turns out not to be sharp.  A
round ball with volume $V(U)$ has k-width roughly $V(U)^{k/n}$. 
If $V(U)^{1/n}$ is much larger than $R_1$, the round ball with
volume $V(U)$ does not fit inside of the rectangle $R$.  What are
the subsets of $R$ that maximize the k-width for a given volume? 
One candidate is a rectangle of dimensions $R_1 \times ... \times
R_l \times S \times ... \times S$, where $l \le k-1$ and $S$ is
between $R_l$ and $R_{l+1}$.  We call this set $U_0$.  The volume
$V(U_0)$ is equal to $R_1 ... R_l S^{n-l}$ and the k-width
$W_k(U_0)$ is approximately $R_1 ... R_l S^{k-l}$.  Solving for
$S$ in terms of $V(U_0)$ and plugging in, we see that $W_k(U_0)$
is roughly equal to $(R_1 ... R_l)^{(n-k)/(n-l)}
V(U_0)^{(k-l)/(n-l)}$.  It turns out that $U_0$ has roughly the
largest k-width among all subsets of $R$ with its volume.  We now
prove that any subset $U$ of $R$ obeys the inequality $W_k(U) <
C(n) (R_1 ... R_l)^{(n-k)/(n-l)} V(U)^{(k-l)/(n-l)}$.  This
inequality becomes roughly an equality when $U = U_0$.

\begin{prop} If $U$ is an open set contained in R, then for
each integer l in the range $0 \le l \le k$, the following
inequality holds.  

$$W_k(U) < C(n) (R_1 ... R_l)^{(n-k)/(n-l)} V(U)^{(k-l)/(n-l)}.$$

\end{prop}

\proof When $l=0$, this inequality reduces to the width-volume
inequality.  When $l=k$, this inequality says that the width of
U is less than $C(n) R_1 ... R_k$.  Since $U$ is a subset of R,
the width of $U$ is at most the width of R, and this inequality
follows.  Now we turn to the intermediate values of l.

Let f be the function on the (n-l)-dimensional rectangle $R_{l+1}
\times ... \times R_n$ with $f(y)$ equal to $(R_1 ... R_l)^{-1}$
times the volume of $U \cap [0, R_1] \times ... \times [0, R_l]
\times \{ y \}$.  In other words, if $U$ contains all of $[0, R_1]
\times ... \times [0, R_l] \times \{ y \}$, then $f(y)$ will be
1, and if $U$ contains half of that region, $f(y)$ will be $1/2$. 

The function f is compactly supported, and $0 \le f \le 1$. 
Applying the width-volume inequality for functions to the
function f, we see that the (k-l)-width of f is bounded by $C(n)
(\int f)^{\frac{k-l}{n-l}}$.  This expression is equal to
$C(n) [(R_1 ... R_l)^{-1} V(U)]^{\frac{k-l}{n-l}}$.  According to
the definition of (k-l)-width, there is a family F
of (k-l)-cycles sweeping out the support of f, so that
the integral of f over each cycle $F(p)$ is bounded by this
expression.  We define a family $F'$ of k-cycles
sweeping out R.  The family $F'$ has the same parameter space as F,
and we define $F'(p) = F(p)
\times [0, R_1] \times ... \times [0, R_l]$.  The volume of U
intersected with a cycle $F'(p)$ is bounded by $(R_1 ... R_l)$
times the integral of f over the corresponding cycle $F(p)$. 
Therefore, the k-width of $U$ is bounded by $C(n) (R_1 ...
R_l)^{\frac{n-k}{n-l}} V(U)^{\frac{k-l}{n-l}}$.
\endproof

This proposition allows us to estimate how many disjoint wide
sets can be packed into a rectangle.  We define $P_{k,N}(U)$ to
be the supremal W so that there exist N disjoint subsets $U_i
\subset U$ each with k-width at least W.  The letter P stands
for packing-width.  For a rectangle $R$, we can estimate
$P_{k,N}(R)$ up to a constant factor C(n).  The formula is a
little complicated, but the geometric meaning is that cutting a
rectangle $R$ into rectangular grids gives roughly the optimal
packings.

\begin{prop} Suppose $R$ is an
n-dimensional rectangle with dimensions $R_1 \le ... \le R_n$. 
Then, up to a factor of $C(n)$, $P_{k,N}(R)$ is equal to the infimum
of the following expression over all integers l in the range $0
\le l \le k$:

$$R_1 ... R_l (R_{l+1} ... R_n)^{\frac{k-l}{n-l}}
N^{-\frac{k-l}{n-l}}.$$
\end{prop}

\proof First we prove that $P_{k,N}(R)$ is bounded above by each
of the expressions in the proposition.  Let $U_i$ be N disjoint
subsets of R.  One of them must have volume at most $R_1 ... R_n
/N$.  Applying the width volume inequality for rectangles, we see
that this set has k-width less than $C(n) (R_1 ...
R_l)^{\frac{n-k}{n-l}} (R_1 ... R_n N^{-1})^{\frac{k-l}{n-l}}$,
for each l between 0 and k.  Expanding this expression, we get
$R_1 ... R_l (R_{l+1} ... R_n)^{\frac{k-l}{n-l}}
N^{-\frac{k-l}{n-l}}$.  This is the inequality we wanted to
prove.

It remains to show that the packing-width is at least as great as
this expression.  To do this we will use the packing formed by
cutting $R$ along a rectangular grid.  First we consider the case
$N \le (R_n ... R_{k+1}) / R_k^{n-k}$.  In this case, we can find
N disjoint subrectangles in $R$ each of dimension $R_1 \times ...
\times R_k \times R_k \times ... \times R_k$.  Each of these
rectangles has k-width roughly $R_1 ... R_k$.  Therefore,
$P_{k,N}(R) \ge c(n) R_1 ... R_k$.  Since the k-width of $R$ is at
most $R_1 ... R_k$, it follows that $P_{k,N}(R) \le R_1 ... R_k$,
and this lower bound is sharp up to a constant factor.  Second we
consider the main case that $N \ge (R_n ... R_{k+1}) /
R_k^{n-k}$.  In this case, we can find N disjoint parallel
rectangles in $R$ each with dimensions $R_1 \times ... \times R_l
\times $S$ \times ... \times S$, for some number $S$ in the range
$R_l \le $S$ \le R_{l+1}$, where $l < k$.  Moreover, we can choose
these rectangles so that they fill up a good portion of the total
volume of R.  In other words, $N R_1 ... R_l S^{n-l} > c(n) R_1
... R_n$.  Solving for S, we see that $S > c(n) (R_n ...
R_{l+1})^{\frac{1}{n-l}} N^{- \frac{1}{n-l}}$.  Now the k-width
of each rectangle is $R_1 ... R_l S^{k-l}$, so we conclude that
$P_{k,N}(R) > c(n) R_1 ... R_l (R_{l+1} ...
R_n)^{\frac{k-l}{n-l}} N^{-\frac{k-l}{n-l}}$.  This lower bound
comes within a constant factor of one of the upper bounds we
proved in the first paragraph.  Therefore, we have determined
$P_{k,N}(R)$ up to a factor $C(n)$. \endproof

\section{Estimates of k-dilation}

In this section, we will estimate the k-dilation of degree 1 maps
between certain domains in Euclidean space, especially
rectangles.  The first estimate follows from our knowledge of the
k-width of rectangles, and the more refined estimates follow from
our knowledge of the packing widths of rectangles.  We begin by
reviewing the definition of k-dilation and some of its basic
properties.

Recall that a piecewise smooth map f has k-dilation at most
$\lambda$ if f maps each k-dimensional submanifold of the domain
with volume V to an image with volume at most $\lambda V$.  The
k-dilation of f can also be expressed in terms of the derivative
$df$.  When $k=1$, the 1-dilation of f is equal to its Lipshitz
constant, which is equal to the supremum of $|df|$.  We now
generalize this result to all values of k.  

If f maps M to N, then $df$ at a point m maps $TM_m$ to
$TN_{f(m)}$.  Taking the k-fold exterior power of this map gives
a map $\Lambda^k df_m$ from $\Lambda^k TM_m$ to $\Lambda^k
TN_{f(m)}$.  By $|\Lambda^k df_m|$, we denote the operator norm
of this linear map.  In other words, this norm is the maximum
over all unit k-vectors v in $\Lambda^k TM_m$ of the norm
$|\Lambda^k df_m (v)|$.

\begin{prop} The k-dilation of a piecewise smooth map f is equal
to the supremum of $|\Lambda_k df_m|$ as m varies over M.
\end{prop}

\proof If the original derivative $df_m$ has singular values $0 \le s_1
\le ... \le s_n$, corresponding to singular vectors $v_1, ...,
v_n$, then the singular values of $\Lambda^k df$ are given by all
products of k distinct numbers $s_i$, and the
singular k-vectors are given by the wedge products of
the corresponding vectors $v_i$.  Therefore, it follows that the norm
$|\Lambda^k df_m|$ is equal to $s_{n-k+1} ... s_n$.  Moreover,
this singular value corresponds to a simple k-vector $v_{n-k+1}
\wedge ... \wedge v_n$.  Taking a small disk near m in the
plane spanned by $v_{n-k+1}, ..., v_n$, we see that the
k-dilation of f is at least $|\Lambda^k df_m|$.  On the other
hand, the linear map $df_m$ stretches the volume of each
k-dimensional disk by at most $|\Lambda^k df_m|$.  By a standard
calculus argument, the k-dilation of f is at most $sup_m
|\Lambda^k df_m|$. \endproof

Formulating the k-dilation in terms of the singular values of the
derivative allows us to show that the k-dilation controls the
$(k+i)$-dilation for all $i > 0$.

\begin{prop} Let f be a piecewise smooth map with k-dilation at
most $\lambda$.  Then for each $i > 0$, the $(k+i)$-dilation of f
is at most $\lambda^{\frac{k+i}{k}}$.
\end{prop}

\proof In the proof of the last proposition, we showed that the
k-dilation of f is equal to the supremum of $s_{n-k+1} ... s_n$. 
Similarly, the $(k+i)$ dilation is equal to the supremum of
$s_{n-k-i+1} ... s_n$.  This expression is bounded by
$s_{n-k+1}^{i} s_{n-k+1} ... s_n$.  Since the k-dilation of f is
at most $\lambda$, $s_{n-k+1} \le \lambda^{1/k}$.  Plugging this
bound into the last expression, we see that the $(k+i)$-dilation
of f is at most $\lambda^{\frac{k+i}{k}}$. \endproof

We now turn to the main problem of this section.  Let $U$ and $V$
be connected bounded open sets in $\mathbb{R}^n$ with piecewise smooth
boundaries.  Let $D_k(U,V)$ denote the infimal k-dilation of a
degree 1 map from the pair $(U, \partial U)$ to the pair $(V,
\partial V)$.  How can we estimate $D_k(U,V)$?

To give some context, we first consider the more familiar cases
when k is equal to n or to 1.  If $k=n$, we can get sharp
estimates by using Moser's theorem for inducing differential
forms.  Suppose that $U$ and $V$ are diffeomorphic.  Let $dvol_U$
denote the volume form on $U$ and $dvol_V$ denote the volume form
on $V$.  Moser proved in \cite{M} that there is a diffeomorphism
$\phi$ from $U$ to $V$ so that $\phi^* dvol_V = \mu dvol_U$, where
$\mu$ is the ratio $\textrm{Volume}(V)/\textrm{Volume}(U)$.  This
diffeomorphism has n-dilation $\mu$.  Since a degree 1 map is
surjective, any degree 1 map from $U$ to $V$ must have n-dilation at
least the ratio $\textrm{Volume}(V)/\textrm{Volume}(U)$.  This
result is very satisfactory, but it has no analogue for $k < n$. 
Next we consider the case $k=1$.  For complicated domains $U$ and
$V$, our problem may be difficult even for $k=1$.  The
distinguishing feature of $k=1$ is that we have a brute force
approach which is not available for higher values of k.  If we
fix bounded open sets $U$ and $V$, then the set of maps from $(U,
\partial U)$ to $(V, \partial V)$ with 1-dilation at most
$\lambda$ is compact in $C^0$.  Therefore, at least in theory,
one can systematically search this class of maps for maps of
degree 1.  This approach can be carried out on a computer if $U$
and $V$ are polyhedra, and it would give an estimate of $D_1(U,V)$
to arbitrary accuracy, although it would be extremely slow.  By
constrast, the set of maps with 2-dilation at most $\lambda$ is
not compact in $C^0$, so that even with unlimited computing time
I don't know how to systematically estimate $D_2(U,V)$ up to a
factor of $10^{100}$.

We will focus on the special case of maps from a rectangle $R$ to
a rectangle $S$.  Even in this special case, the problem is much
harder than I initially expected.  We use the convention that the
rectangle $R$ has dimensions $R_1 \le ... \le R_n$ and the
rectangle $S$ has dimensions $S_1 \le ... \le S_n$.  To make the
algebra simpler, we let $Q_i$ denote the quotient $S_i/R_i$.  We
now prove some lower bounds for k-dilation.

\begin{prop} Suppose that $U$ is a subset of R.  Then $D_k(U,S)$ is
at least $c(n) Q_1 ... Q_k$.
\end{prop}

\proof Since $U$ is a subset of R, $W_k(R) \ge W_k(U)$.  Now if
there is a degree non-zero map from $U$ to $S$ with k-dilation
$\lambda$, then $\lambda W_k(U) \ge W_k(S)$.  Therefore, $\lambda
\ge W_k(S) / W_k(R)$.  But according to Proposition 1.2,
$W_k(S)$ is at least $c(n) S_1 ... S_k$, and $W_k(R)$ is at most $R_1
... R_k$.  Therefore, $\lambda$ is at least $c(n) Q_1 ... Q_k$.
\endproof

We can get more complicated bounds by considering the
packing-widths of $R$ and $S$.

\begin{prop} Suppose that $U$ is a subset of $R$.  Then, for each
integer l from 0 to k, $D_k(U,S)$ is at least $c(n) Q_1 ... Q_l
(Q_{l+1} ... Q_n)^{\frac{k-l}{n-l}}$.
\end{prop}

\proof Since $U$ is a subset of $R$, we have $P_{k,N}(R) \ge
P_{k,N}(U)$ for every k and N.  If there is a map $f$ of non-zero
degree from $U$ to $S$ with k-dilation $\lambda$, then $\lambda
P_{k,N}(U) \ge P_{k,N}(S)$.  To see this, let $S_i$ be N disjoint
subsets of S, each with width at least $P_{k,N}(S) -
\epsilon$.  Then let $U_i$ be the inverse image $F^{-1}(S_i)$. 
The map $f$ restricts to a degree non-zero map from $(U_i, \partial
U_i)$ to $(S_i, \partial S_i)$.  Therefore the k-width of $U_i$
is at least $\lambda^{-1} (P_{k,N}(S) - \epsilon)$.  Since the
sets $U_i$ are disjoint, $\lambda P_{k,N}(U) \ge P_{k,N}(S)$. 
This estimate gives us a lower bound $\lambda \ge P_{k,N}(S) /
P_{k,N}(R)$, for every natural number N.

The value of $P_{k,N}(R)$ is estimated in Proposition 2.3.  Up
to a constant factor C(n), it is equal to the infimum of $R_1 ...
R_l (R_{l+1} ... R_n)^{\frac{k-l}{n-l}} N^{-\frac{k-l}{n-l}}$,
where l lies in the range $0 \le l \le k$.  In particular, we can
consider the case that $N = S_n ... S_{l+1} / S_l^{n-l}$.  Since
there are roughly N disjoint rectangles in $S$ with dimensions $S_1
\times ... \times S_l \times S_l \times ... \times S_l$, the
packing-width $P_{k,N}(S)$ is at least $c(n) S_1 ... S_l
S_l^{k-l}$.  On the other hand, $P_{k,N}(R)$ is at most $C(n) R_1
... R_l (R_{l+1} ... R_n)^{\frac{k-l}{n-l}} [S_n ... S_{l+1} /
S_l^{n-l}]^{-\frac{k-l}{n-l}}$.  Therefore,
$P_{k,N}(S)/P_{k,N}(R)$ is at least $c(n) Q_1 ... Q_l (Q_{l+1} ...
Q_n)^{\frac{k-l}{n-l}}$.

This finishes the proof of the proposition.  The reader can check
that the supremum over N of the quotient
$P_{k,N}(S)/P_{k,N}(R)$ is approximately equal to the maximum
of $Q_1 ... Q_l (Q_{l+1} ... Q_n)^{\frac{k-l}{n-l}}$ for l in the
range $0 \le l \le k$.  Therefore, the packing-width does not
give any further lower bounds for $D_k(R,S)$.  \endproof

The second theorem of this paper is an estimate for $D_{n-1}(R,S)$.

\begin{theorem} Let $R$ and $S$ be n-dimensional rectangles. 
Suppose that $R$ has dimensions $R_1 \le ... \le R_n$ and that S
has dimensions $S_1 \le ... \le S_n$.  Up to a constant factor
C(n), $D_{n-1}(R,S)$ is equal to the supremum of the following
quantities:

$$Q_1 ... Q_l (Q_{l+1} ... Q_n)^{\frac{n-l-1}{n-l}} \eqno{(1)}$$

$$Q_2 ... Q_n. \eqno{(2)}$$

\noindent In equation (1), the number l is allowed to take any
value in the range $1 \le l \le n-1$.
\end{theorem}

For example, if n is 3, then $D_2(R,S)$ is roughly the supremum
of $Q_1 Q_2^{1/2} Q_3^{1/2}$, $Q_1 Q_2$, and $Q_2 Q_3$.

We have already proven the lower bounds in equation 1.  They are
exactly the lower bounds in Proposition 3.4 in case $k=n-1$.

The lower bound in equation 2 is simple.  Suppose that $f$ is a
degree 1 map from $R$ to $S$ with $(n-1)$-dilation $\lambda$.  The
map f restricts to a degree 1 map from the boundary of $R$ to the
boundary of S.  Since this map must be surjective, it follows
that $\lambda \textrm{Volume}(\partial R) \ge
\textrm{Volume}(\partial S)$.  But the volume of the boundary of R
is at most $2n R_2 ... R_n$, and the volume of the boundary of S
is at least $2 S_2 ... S_n$.  Therefore, $\lambda$ is at least
$(1/n) Q_2 ... Q_n$.  This finishes the proof of the lower bounds
on $D_{n-1}(R,S)$.

In the next section, we will construct degree 1 maps showing that
these lower bounds are sharp up to a constant factor.

\section{Maps with small $(n-1)$-dilation}

In this section, we construct a degree 1 map between rectangles
with surprisingly small k-dilation.  After the construction, we
check that the k-dilation of this map can be smaller than the
k-dilation of any linear diffeomorphism by an arbitrarily large
factor.  Then we will finish the proof of Theorem 2, determining
$D_{n-1}(R,S)$ up to a constant factor.

\begin{construction} (The snake map) Let $R$ and $S$ be n-dimensional
rectangles.  Suppose that $n \ge 3$ and that $k$ lies in the
range $2 \le k \le n-1$.  Suppose that $R_i = S_i$ for $i \le
n-k$.  Suppose that $R_{n-k+1} ... R_{n-k+b} \ge S_{n-k+1} ...
S_{n-k+b}$ for every b in the range $1 \le b \le k$.  Then there
is a degree 1 map from $R$ to $S$ with k-dilation less than
$C(n)$.
\end{construction}

\proof We write $R$ as the product $R' \times R''$, where $R' = [0,
R_1] \times ... \times [0, R_{n-k}]$ and $R'' = [0, R_{n-k+1}]
\times ... \times [0, R_n]$.  Similarly, we write $S = S' \times
S''$.  By assumption $R'$ is congruent to $S'$.

Because of the inequalities $R_{n-k+1} ... R_{n-k+b} \ge
S_{n-k+1} ... S_{n-k+b}$, there is a smooth bilipschitz
embedding of $S''$ into $R''$ with bilipschitz constant at
most $C(n)$.  We will need a little bit of room later, so we let $I$
be a smooth bilipschitz embedding of $3 S''$ into $R''$, with
quasi-isometric constant $C(n)$.  (By $3 S''$, we mean the rectangle
$S''$ dilated by a factor of 3 around its center.)  We let A be
the image of I in $R''$.  (This set A corresponds to the set A in the
description of the snake map before this proof.)  Let H
be a smooth function on $3 S''$ which is equal to 1 on the
central $S''$ and is equal to zero on a neighborhood of the
boundary of $S''$.  We can choose H with Lipschitz constant as
close as we like to $S_{n-k+1}^{-1}$.

The function $H \circ I^{-1}$ is defined on the image of I in
$R''$, and it is equal to zero on the boundary of this image.  We
extend this function to all of $R''$ by setting it equal to zero
on the complement of the image of I.  We call the resulting
function $\bar H$.  We denote a point in $R$ by $(x', x'')$, where
$x'$ lies in $R'$ and $x''$ lies in $R''$.  We define $\Phi_1(x',
x'') = (\bar H(x'') x', x'')$.  If we differentiate $\Phi_1$, we
find that the norm of the derivative is bounded by the sum $\sup
|H(x'')| + \sup |x'| sup |\nabla \bar H|$.  The first of these
expressions is bounded by 1, and the second by $C(n)
R_{n-k}/S_{n-k+1}$.  Because of our assumptions about the
dimensions of $R$ and S, we have $R_{n-k} = S_{n-k} \le S_{n-k+1}$,
and so the map $\Phi_1$ has Lipschitz constant less than C.  The
image $\Phi_1(R)$ is contained in $R' \times A \cup \{ 0 \}
\times R''$.  We call this set Q.

The next step of our construction is to retract the region Q onto
$R' \times A$.  To do this, we first pick a retraction $\phi_2$
from $R''$ to the image of I.  We choose $\phi_2$ so that it maps
the complement of A onto the boundary of A.  We also assume that
$\phi_2$ is piecewise smooth.  Next, we define $\Phi_2(x',x'') =
(x', \phi_2(x''))$.

The map $\Phi_2$ has large k-dilation on R, but its restriction
to Q has k-dilation 1.  On the intersection of Q with the region
$R' \times A$, $\Phi_2$ is the identity, and so it has k-dilation
1.  The complement of this region in Q is given by the conditions
$x' = 0$ and $x'' \in A^c$, where $A^c$ denotes the complement of
A in $R''$.  The map $\Phi_2$ takes this $k$-dimensional region
into the $(k-1)$-dimension region given by the conditions $x'=0$,
and $x'' \in \partial A$.  Therefore, $\Phi_2$ has k-dilation
zero on the second part of Q.  All together, the map $\Phi_2$ has
k-dilation 1.

Next, we define a map $\Phi_3$ from the region $R' \times A$ to
$S' \times 3 S''$.  This map is defined by $\Phi_3(x', x'') =(x',
I^{-1}(x''))$.  It has Lipschitz constant at most $C(n)$.  The
composition $\Phi_3 \circ \Phi_2 \circ \Phi_1$ is a map of $R$
into $S' \times 3 S''$, with k-dilation less than $C(n)$.  The
rectangle $S' \times 3 S''$ contains the rectangle $S = S' \times
S''$.  Since $S$ is convex, there is a retraction $\Phi_4$ from
$S' \times 3S''$ to $S$ with Lipschitz constant 1.

The composition $\Phi_4 \circ \Phi_3 \circ \Phi_2 \circ \Phi_1$
is a degree 1 map from $(R, \partial R)$ to $(S, \partial S)$
with k-dilation less than $C(n)$.   \endproof

Remark: With a little more work, it is possible to construct a PL
isomorphism from $R$ to $S$ with k-dilation less than $C(n)$.

We now give an example to show that the snake map badly
outperforms the linear map for some rectangles.  Let $Lin_k(R,S)$
denote the smallest k-dilation of any linear diffeomorphism from
$R$ to $S$.

\begin{prop} For each $n \ge 3$ and each k in the range $2 \le k
\le n-1$, there are n-dimensional rectangles $R$ and $S$ which
make the ratio $Lin_k(R,S) / D_k(R,S)$ arbitrarily large.
\end{prop}

\proof Let $R$ be the rectangle with dimensions $R_1 = ... =
R_{n-2} = \epsilon$, and $R_{n-1} = R_n = 1$.  Let $S$ be the
rectangle with dimensions $S_1 = ... = S_{n-1} = \epsilon$ and
$S_n = \epsilon^{-1}$.  Construction 2 gives a degree 1 map from
R to $S$ with 2-dilation at most $C(n)$.  Next we give a lower
bound for $Lin_{n-1}(R,S)$.  Any linear diffeomorphism from $R$ to
S takes each hyperface of $R$ onto a hyperface of S.  The rectangle
R has 4 hyperfaces with volume $\epsilon^{n-2}$.  On the other
hand, the rectangle $S$ has $2(n-1)$ hyperfaces with volume
$\epsilon^{n-3}$ and only 2 hyperfaces with volume
$\epsilon^{n-2}$.  Any linear diffeomorphism from $R$ onto $S$ must
take a face of $R$ with volume $\epsilon^{n-2}$ onto a face of S
with volume $\epsilon^{n-3}$.  Therefore, it must have
(n-1)-dilation at least $\epsilon^{-1}$, and we conclude that
$Lin_{n-1}(R,S) \ge \epsilon^{-1}$.

According to Proposition 3.2, a map with k-dilation $\lambda$ has
$(k+i)$-dilation at most $\lambda^{\frac{k+i}{k}}$.  Therefore,
$D_k(R,S)$ is less than $C(n)$ for all $k \ge 2$.  By the same
argument, $Lin_k(R,S)$ is at least $\epsilon^{-\frac{k}{n-1}}$
for all $k \le n-1$.  Combining these estimates, we see that
$Lin_k(R,S) / D_k(R,S)$ may be arbitrarily large for all k in the
range $2 \le k \le n-1$. \endproof

Using the snake map, we now finish the proof of Theorem 2.

\proof By composing snake
maps and linear maps, we will construct enough degree 1 maps to
prove the theorem.  We begin with the case $n=3$. 

By scaling the rectangle S, we can assume that the lower bound
for $D_2(R,S)$ given in Theorem 2 is equal to 1.  In other words,
we suppose $R_1 R_2 > S_1 S_2$, $R_1^2 R_2 R_3 > S_1^2 S_2 S_3$,
and $R_2 R_3 > S_2 S_3$.  Under these assumptions, we need to
construct a degree 1 map from $R$ to $S$ with 2-dilation less than C. 
We do so in three cases.

If $R_1 < S_1$, then we define a 2-contracting linear
diffeomorphism from $R$ to T, with $T_1 = S_1$, $T_2 = R_2
R_1/S_1$, and $T_3 = R_3 R_1/ S_1$.  (The length $T_2$ is indeed
bigger than $T_1$ because $R_1 R_2 > S_1 S_2$.)  Using the
first two equations in the list above, we see that $T_2 > S_2$
and $T_2 T_3 > S_2 S_3$.  Therefore, there is a snake map from T
to $S$ with 2-dilation less than C.

If $R_1 \ge S_1$ but $R_2 < S_2$, then we define a 2-contracting
linear diffeomorphism from $R$ to T, with $T_1 = R_1 R_2/ S_2$,
$T_2 = S_2$, and $T_3 = R_3 R_2/ S_2$.  (The length $T_3$ is
indeed bigger than $T_2$ because $R_2 R_3 > S_2 S_3$.)  Since
$R_1 R_2 > S_1 S_2$, $T_1 > S_1$.  Since $R_2 R_3 > S_2 S_3$,
$T_3 > S_3$.  Therefore, there is a 1-contracting linear
diffeomorphism from T to S.

If $R_1 \ge S_1$ and $R_2 \ge S_2$, since we have assumed that
$R_2 R_3 \ge S_2 S_3$, there is a snake map from $R$ to $S$ with
2-dilation less than C. 

We now turn to the case of higher dimensions.  As in the
three-dimensional case, we can scale $S$ so that the lower bound on
$D_{n-1}(R,S)$ is equal to 1.  In other words, we can assume that
the rectangles $R$ and $S$ obey the following list of inequalities
denoted $(*)$:

$$R_1 ... R_l (R_{l+1} ... R_n)^{\frac{n-l-1}{n-l}} > S_1 ... S_l
(S_{l+1} ... S_n)^{\frac{n-l-1}{n-l}}. \eqno(*1)$$

$$R_2 ... R_n \ge S_2 ... S_n. \eqno(*2)$$

\noindent Equation $(*1)$ holds for every integer l in the range
$1 \le l \le n-1$.

Assuming $(*)$, we have to construct a degree 1 map from $R$ to
$S$ with $(n-1)$-dilation at most $C(n)$.  The maps we will
construct will have the following structure.  First, there will
be a snake map from $R$ to an intermediate rectangle $T$, with
(n-1)-dilation at most $C(n)$.  Then there will be an
(n-1)-contracting linear diffeomorphism from $T$ to $S$. 
Choosing the rectangle $T$ and constructing the two maps is very
tedious, but it requires only elementary algebra.

We begin with the special case $R_1 = S_1$.  Because of the messy
algebra, from now on we use $C$ to denote a constant that depends
only on n.

\begin{lemma} Suppose that $R$ and $S$ are n-dimensional rectangles
obeying $(*)$.  Also suppose that $R_1 = S_1$, and suppose that
Theorem 2 holds for rectangles of dimension n-1.  Then there is a
degree 1 map from $R$ to $S$ with (n-1)-dilation at most C.
\end{lemma}

\proof Since $R_1 = S_1$, the inequalities in $(*)$ imply that
for each l in the range $1 \le l \le n-1$,  
$R_2 ... R_l (R_{l+1} ... R_n)^{\frac{n-l-1}{n-l}} \ge S_2
... S_l (S_{l+1} ... S_n)^{\frac{n-l-1}{n-l}}$.

We are going to prove something a little more general in order to
do an inductive argument.  For each p in the range $1 \le p \le
n-1$, let $C(p)$ denote the following list of conditions.

1. $R_1 = S_1$.

2. $R_2 ... R_a \ge S_2 ... S_a$ for a in the range $2 \le
a \le p$,

3. $R_2 ... R_l (R_{l+1} ... R_n)^{\frac{n-l-1}{n-l}} \ge S_2 ...
S_l (S_{l+1} ... S_n)^{\frac{n-l-1}{n-l}}$ for l in the range $p
\le l \le n-1$, and

4. $R_2 ... R_n \ge S_2 ... S_n$.

Our hypotheses are exactly $C(1)$.  We are going to prove that
for every p, $C(p)$ implies that there is a degree 1 map from R
to $S$ with $(n-1)$-dilation at most C.  The point of introducing
all of these new conditions is that we can make an inductive
argument, starting with $C(n-1)$ and working our way down to $C(1)$.

To anchor the induction we prove that $C(n-1)$ implies that there
is a degree 1 map from $R$ to $S$ with $(n-1)$-dilation at most C. 
The condition $C(n-1)$ says exactly that $R_1 = S_1$, and that
$R_2 ... R_a \ge S_2 ... S_a$ for every a.  Under these
conditions, our construction gives a snake map from $R$ to $S$ with
$(n-1)$-dilation at most C.

We now turn to the inductive step.  We may assume that $C(q)$
implies a good map for every $q$ greater than p, and we have to
prove that $C(p)$ implies the existence of a good map.

Suppose that $R$ and $S$ satisfy $C(p)$.  Let b be the smallest
number for which $R_2 ... R_b < S_2 ... S_b$.  If there is no
such b, then $R$ and $S$ actually satisfy $C(n-1)$, and so there
is a snake map from $R$ to $S$ with (n-1)-dilation at most $C$. 
Because of condition 4 we know that b is not equal to n, and
because of the condition 3 with $l=n-1$, we know that b is not
equal to n-1.  On the other hand, $b$ must be greater than $p$. 
Therefore, b lies in the range $p < b \le n-2$.

We begin with case $b=2$, which is a bit easier than the general
case.  If $b=2$, then we must have had $p=1$.  There is an
(n-1)-contracting linear diffeomorphism from $T$ to $S$, where
$T_1 = S_1$, $T_2 = R_2$, and $T_i = S_i (S_2/R_2)^{1/(n-3)}$ for
all $i \ge 3$.  (The length $T_2$ is at least $T_1$ because $T_1
= S_1 = R_1 \le R_2 = T_2$.)  We check that $R$ and $T$ obey
$C(2)$.  Condition 1 follows because $R_1 = S_1 = T_1$. 
Condition 2 follows because $R_2 = T_2$.  A computation shows
that for $l$ in the range $2 \le l \le n-1$, $T_2 ... T_l
(T_{l+1} ... T_n)^{\frac{n-l-1}{n-l}} = S_2 ... S_l (S_{l+1} ...
S_n)^{\frac{n-l-1}{n-l}}$.  Therefore, condition 3 of $C(2)$
holds.  Finally, condition 4 follows from condition 3 for l=2
along with the equality $T_2 = R_2$.  Since $p=1$, our inductive
hypothesis is that $C(2)$ implies the existence of a good map. 
Therefore we may conclude there is a good map from $R$ to
$T$.  Composing with the linear map from $T$ to $S$ gives a
good map from $R$ to $S$.

Now we deal with the more general case that $b > 2$.  In this
case, we apply an (n-1)-expanding linear transformation to $S$
that leaves $S_1$ through $S_{b-1}$ invariant, decreases $S_b$ by
some factor $\lambda$, and increases $S_{b+1}$ through $S_n$ by a
factor $\lambda^{\frac{1}{n-b-1}}$.  (At this stage, we use the
fact that $b < n-1$.)  We choose $\lambda$ so that the image
rectangle $S'$ has either $R_2 ... R_b = S'_2 ... S'_b$, or
$S'_{b-1} = S'_b$, whichever requires a smaller value of
$\lambda$.  This choice of $\lambda$ ensures that the dimensions
of $S'$ are still labelled in increasing order.  If $R_2 ... R_b
= S'_2 ... S'_b$, then we can stop, but if $S'_{b-1} = S'_b$, we
have to proceed with another (n-1)-expanding linear
transformation.  In this case, we then apply an (n-1)-expanding
linear transformation to $S'$ that leaves $S'_1$ through
$S'_{b-2}$ invariant, decreases $S'_{b-1}$ and $S'_b$ equally,
and increases all the other directions equally.  We choose the
factor of stretching so that the image rectangle $S''$ has either
$R_2 ... R_b = S''_2 ... S''_b$ or $S''_{b-2} = S''_{b-1} =
S''_b$.  In the latter case, we then apply a linear
transformation that decreases $S''_{b-2}, S''_{b-1},$ and
$S''_b$, and so on.  Because $R_2 > S_2$, this process
terminates.  We call the final rectangle in this chain of
(n-1)-expanding diffeomorphisms $T$.

At the end of the process, we have an equality $R_2 ... R_b = T_2
... T_b$.  We check that the rectangles $R$ and $T$ obey the
condition $C(b)$.  Condition 1 follows because $R_1 = S_1 = T_1$. 
Condition 2 follows because $R_2 ... R_b = T_2 ... T_b$, and
because for each $a$ less than $b$, $T_2 ... T_a \le S_2 ... S_a
\le R_2 ... R_a$.  A calculation shows that $T_2 ... T_l
(T_{l+1} ... T_n)^{\frac{n-l-1}{n-l}} = S_2 ... S_l (S_{l+1}
... S_n)^{\frac{n-l-1}{n-l}}$ for $l \ge b$.  Therefore,
condition 3 of $C(b)$ holds for $R$ and $T$.  Finally condition 4
follows from the case $l=b$ of condition 3 along with the equality
$R_2 ... R_b = T_2 ... T_b$.  Since $b$ is greater than $p$, our
inductive hypothesis tells us that there is a degree 1 map from
$R$ to $T$ with (n-1)-dilation at most C.  Since there is an
(n-1)-contracting linear diffeomorphism from $T$ to $S$, we can
compose them to get a degree 1 map from $R$ to $S$ with
$(n-1)$-dilation at most $C$.

This argument proves the inductive step and hence the lemma.
\endproof

With the help of the lemma, we can now prove Theorem 2 for all
pairs of rectangles $R$ and $S$.  We assume that $R$ and $S$ obey $(*)$,
and we need to construct a degree 1 map from $R$ to $S$ with
(n-1)-dilation at most C.

If $R_1 < S_1$, then we proceed as follows.  There is an
(n-1)-contracting linear map from $T$ to $S$, where $T$ has dimensions
$T_1 = R_1$, and $T_i = S_i (S_1/R_1)^{\frac{1}{n-2}}$ for all
other i.  A calculation shows that $T_1 ... T_l (T_{l+1}
... T_n)^{\frac{n-l-1}{n-l}} = S_1 ... S_l (S_{l+1} ...
S_n)^{\frac{n-l-1}{n-l}}$ for each l in the range $1 \le l \le
n-1$.  Since $R$ and $S$ obey $(*)$, it follows that $R_1 ... R_l
(R_{l+1} ... R_n)^{\frac{n-l-1}{n-l}}
\ge T_1 ... T_l (T_{l+1} ... T_n)^{\frac{n-l-1}{n-l}}$ for each l
in the range $1 \le l \le n-1$.  Since $R_1 = T_1$, the last
equation in the case $l=1$ implies that $R_2 ... R_n \ge T_2 ...
T_n$.  In other words, $R$ and $T$ obey $(*)$.  Since $R_1 = T_1$, we
can apply Lemma 4.1, which tells us that there is a degree 1 map
from $R$ to $T$ with $(n-1)$-dilation less than C.  Composing with the
linear map from $T$ to $S$ gives a degree 1 map from $R$ to $S$ with
$(n-1)$-dilation less than C, which is what we wanted to prove.

If $R_1 > S_1$, then we proceed as follows.  We apply an
(n-1)-expanding linear transformation to $S$ which decreases $S_2$
and increases all other directions of $S$ equally until either $S_1
= R_1$ or $S_1 = S_2$.  In the second case, we apply an
(n-1)-expanding linear transformation to $S$ which decreases $S_3$
and increases all other directions of $S$ equally until either $S_1
= R_1$ or $S_1 = S_2 = S_3$, and so on.  We continue this process
until either $S_1 = R_1$ or $S_1 = S_n < R_1$.  In the latter
case, there is a contracting linear diffeomorphism from $R$ to $S$.

In the former case, we call the final rectangle in this chain of
(n-1)-expanding linear diffeomorphisms $T$.  We have to check
that $R$ and $T$ obey $(*)$.  Suppose that the last diffeomorphism
was decreasing $S_{b+1}$ and increasing $S_i$ for all other i. 
We have $R_1 = T_1 = T_2 = ... = T_b$.  If $l$ is at least $b+1$,
then a calculation shows that $T_1 ... T_l (T_{l+1} ...
T_n)^{\frac{n-l-1}{n-l}} = S_1 ... S_l (S_{l+1} ...
S_n)^{\frac{n-l-1}{n-l}}$.  This equation implies that $R_1 ...
R_l (R_{l+1} ... R_n)^{\frac{n-l-1}{n-l}} \ge T_1 ... T_l
(T_{l+1} ... T_n)^{\frac{n-l-1}{n-l}}$ for $l \ge b+1$. 
Another short calculation shows that $T_2 ... T_n = S_2 ... S_n$. 
Therefore, $R_2 ... R_n \ge T_2 ... T_n$.  Since $T_1 = T_b =
R_1$, it follows that $R_1 ... R_l \ge T_1 .. T_l$ for all l in
the range $1 \le l \le b$.  Since $R_1 = T_1$ and $R_2 ... R_n
\ge T_2 ... T_n$, it follows that $R_1 ... R_n \ge T_1 ... T_n$. 
Combining the last two inequalities, it follows that $R_1 ... R_l
(R_{l+1} ... R_n)^{\frac{n-l-1}{n-l}} \ge T_1 ... T_l (T_{l+1}
... T_n)^{\frac{n-l-1}{n-l}}$, for all $b$ in the range $1 \le l
\le b$.  Assembling all these inequalities, we see that $R$ and
$T$ obey $(*)$.

Since $R_1 = T_1$, we can apply Lemma 4.1, which tells us that
there is a degree 1 map from $R$ to $T$ with $(n-1)$-dilation
less than C.  Composing with the linear diffeomorphism from $T$
to $S$ gives a degree 1 map from $R$ to $S$ with $(n-1)$-dilation
less than C.  \endproof

\section{Appendix: Dividing area and (n-1)-width}

In this appendix, we briefly consider estimating the width of a
Riemannian manifold $(M^n, g)$.  The main result is that for each
$n \ge 2$ a closed oriented Riemannian n-manifold with volume 1
may have arbitrarily large (n-1)-width.  Therefore, the
width-volume inequality does not extend to Riemannian manifolds. 

Our results are based on isoperimetric inequalities.  Let $(M^n,
g)$ be a closed Riemannian manifold.  We define the dividing area
of $(M, g)$ to be the infimum of the volume of $\partial A$ as
$A \subset M$ varies over all open sets with volume between $(1/4)
\textrm{Volume}(M)$ and $(3/4) \textrm{Volume}(M)$.  Any
isoperimetric inequality on $(M,g)$ will lead to a lower bound for
its dividing area.

\begin{prop} Let $(M,g)$ be a closed oriented Riemannian
n-manifold with dividing area $A(M,g)$.

$$W_{n-1}(M,g) \ge (1/2) A(M,g).$$
\end{prop}

Recall that the Cheeger isoperimetric constant $h(M,g)$ is equal
to the supremum of $\textrm{vol}(A) / \textrm{vol} (\partial A)$
as A ranges over all open sets in $(M,g)$ with volume at most
$(1/2) \textrm{vol}(M)$.  In \cite{B}, Brooks constructed
examples of Riemannian manifolds with arbitrarily large volume
and with Cheeger constant bounded below.  We state Brooks's
result as a theorem.

\begin{reftheorem} (Brooks) For each $n \ge 2$, there is a
sequence of closed oriented Riemannian n-manifolds $(M_i, g_i)$
with volume $V_i$ tending to infinity and $h(M_i, g_i) > 1$. 
\end{reftheorem}

The Riemannian manifolds $(M_i, g_i)$ are finite coverings of a
fixed base manifold with an interesting fundamental group.  The
base manifold could be hyperbolic, or a higher-rank symmetric
space.  Because of the lower bound on the Cheeger constant, the
dividing area of $(M_i, g_i)$ must be at least $(1/4) V_i$. 
According to Proposition 5.1, the (n-1)-width of $(M_i, g_i)$ is
at least $(1/8) V_i$.  Now we rescale $(M_i, g_i)$ so that it has
volume 1.  The rescaled version has (n-1)-width at least $(1/8)
V_i^{1/n}$ which tends to infinity.

In Brooks's examples, the topology of $M_i$ is unbounded. There
are other examples on the n-sphere for $n \ge 3$.  These examples
follow from the isoperimetric estimates of Burago and Ivanov,
proven in \cite{BI}.  A small modification of their construction
gives the following theorem.

\begin{reftheorem} (Burago, Ivanov) For each $n \ge 3$ and each
$\epsilon > 0$, there is a metric g on $S^n$ which obeys the
following isoperimetric inequality.  For any open set $A \subset
S^n$ with volume between $(1/4) \textrm{Volume}(S^n, g)$ and
$(3/4) \textrm{Volume}(S^n, g)$,

$$\textrm{Volume}(A)^{\frac{n-1}{n}} < \epsilon
\textrm{Volume}(\partial A).$$

\end{reftheorem}

Since this isoperimetric inequality is scale invariant, we can
also scale the metrics so that they all have volume 1.  In this
case, we have $(1/4) \le \textrm{Volume}(A)^{\frac{n-1}{n}}
\le \epsilon \textrm{Volume}(\partial A)$.  Therefore, these
metrics have dividing area at least $(1/4) \epsilon^{-1}$.  By
Proposition 5.1, we conclude that $W_{n-1}(S^n, g) \ge (1/8)
\epsilon^{-1}$.

Now we turn to the proof of Proposition 5.1.

\proof Let F be a family of (n-1)-cycles sweeping out $(M,g)$
with width almost $W_{n-1}(M,g)$.  Without loss of generality we
may assume that the parameter space of F is a circle.  As in
section 1, we perturb F to get a complex of cycles, where each
0-simplex gets mapped to an (n-1)-cycle in M with mass at most
$W_{n-1}(M,g) + \delta$ and each 1-simplex gets mapped to an
n-chain with mass at most $\delta$.

Let $I$ be an interval of the circle consisting of a union of
1-simplices, and consider the union of the corresponding
n-chains, which we denote $C(I)$.  We consider the image of
$C(I)$, which is a set in $(M,g)$.  If $I$ is a single 1-simplex,
this set has volume at most $\delta$.  If $I$ is the whole
circle, then this set is all of $(M,g)$.  Since adding a
1-simplex only slightly changes the volume of this set, we can
find an interval $I$ so that the volume of the image of $C(I)$ is
close to $(1/2) \textrm{Volume}(M,g)$.  Therefore, the volume of
the boundary of this image is at least the dividing area of
$(M,g)$.  Now let $v_1$ and $v_2$ be the boundary vertices of I. 
The boundary of the image of $C(I)$ is contained in the union
$C(v_1) \cup C(v_2)$.  Therefore, one of these two cycles must
have mass at least $(1/2) A(M,g)$.  Taking $\delta \rightarrow 0$
finishes the proof.  \endproof

\section{Appendix: Falconer's estimate for the linear k-width}

In \cite{Fal}, Falconer proved the following theorem, which we have
reformulated in our language.

\begin{reftheorem} (Falconer) Let $U$ be a bounded open set in
$\mathbb{R}^n$.  Suppose that $k > n/2$.  Then there is a family
of parallel k-planes, each interesecting $U$ in a region of
k-volume at most $C(n) \textrm{Volume}(U)^{k/n}$.
\end{reftheorem}

The proof is based on Fourier analysis.  We give a sketch in the
simplest case: $k=2$, $n=3$.  By a scaling argument, we can
assume that the volume of $U$ is 1, and we let $f$ denote the
characteristic function of $U$.  Then we consider the Fourier
transform of $f$.  Because $\| f \|_1 = 1$, we have $\| \hat f
\|_\infty \le 1$.   Because $\| f \|_2 = 1$, the Plancherel
theoerem tells us that $\| \hat f \|_2 = 1$.  We write this last
equation in polar coordinates.

$$\int_{S^2} \int_0^{\infty} |\hat f(\theta, r)|^2 r^2 dr d\theta = 1.$$

It's convenient to expand the polar coordinates so that the
radius takes values on the whole real line by identifying
$(\theta, -r)$ with $(-\theta, r)$.

$$\int_{S^2} \int_{- \infty}^\infty |\hat f(\theta, r)|^2 r^2 dr d\theta =
2.$$

Since the unit sphere has area $4 \pi$, we conclude that for some
choice of $\theta$, we have the following inequality.

$$\int_{-\infty}^{\infty} |\hat f(\theta, r)|^2 r^2 dr \le
\frac{1}{2 \pi}. \eqno(*)$$

Now the main idea of the proof is that for fixed $\theta$, the
function $\hat f(\theta, r)$ encodes the integrals of f over all
the planes perpendicular to $\theta$.  This idea appears in the
theory of the Radon transform.  The problem is rotationally
invariant, so we may assume that $\theta = (0, 0, 1)$.  Define
the averaged function $F(z)$ to be $\int_{\mathbb{R}^2} f(x,y,z)
dx dy$.  Then an elementary calculation shows that $\hat
F(\xi) = \hat f(\theta, \xi)$.  Using equation $(*)$, we
can estimate $\hat F$ sufficiently well to bound $\sup_z |F(z)|$.

In addition to $(*)$, we also know that $|\hat F(\xi)| = |\hat
f(\theta, \xi)| \le 1$ everywhere.  We combine these
inequalities.

$$\int_{-\infty}^{\infty} |\hat F(\xi)|^2 (1 + |\xi|)^2 d\xi < 5.$$

Next we use the Cauchy-Schwarz inequality to bound
$\int_{- \infty}^{\infty} |\hat F(\xi)| d \xi$.

$$ \int_{-
\infty}^{\infty} |\hat F(\xi)| d \xi = \int_{-
\infty}^{\infty} [|\hat F(\xi)| (1 + |\xi|)] (1 + |\xi|)^{-1}
d\xi \le $$

$$\le [\int _{-\infty}^{\infty} |\hat F(\xi)|^2 (1 + |\xi|)^2
d\xi]^{1/2} [\int_{-\infty}^{\infty} (1 + |\xi|)^{-2}]^{1/2} <
\sqrt{10}.$$

Finally, by the Fourier inversion theorem, we conclude that $\| F
\|_\infty \le \| \hat F \|_1 < \sqrt{10}$.  In other words,
every integral $\int_{\mathbb{R}^2} f(x,y,z) dx dy$ is less than
$\sqrt{10}$.  Since $f$ is the characteristic function of $U$,
the intersection of $U$ with each plane $z = constant$ has area
less than $\sqrt{10}$.

In the general case $k > n/2$, the proof is only slightly more
complicated.  Instead of polar coordinates, one has to average
over the Grassman manifold of k-planes in $\mathbb{R}^n$, and
instead of the Cauchy-Schwarz inequality, one has to use the
Holder inequality.

Falconer's theorem does not extend to the case $k=1$ because of
the Besicovitch example.  I don't know whether it extends to $k$
in the intermediate range $2 \le k \le n/2$.  For more
information, consult \cite{Bo}, \cite{Wol}, and the references
therein.

\end{document}